\documentclass[11pt,a4paper]{article}
\usepackage{fancyhdr}
\usepackage[centertags]{amsmath}
\usepackage{amssymb}
\usepackage{mathrsfs}
\usepackage{mathtools}
\usepackage{hyperref}
\usepackage{array}
\usepackage{breqn,caption,subcaption}

\newtheorem{defn}{Definition}[section]
\newtheorem{ex}{Example}[section]
\begin{document}

\begin{center}

\noindent
{\bf \huge Exact solutions of fractional partial differential equations by Sumudu transform iterative method}

{\Large Manoj Kumar\\}
\textit{Department of Mathematics, Savitribai Phule Pune University,\\ Pune - 411007, India\\
	National Defence Academy, Khadakwasala Pune- 411023, India\\
	mkmath14@gmail.com}
\vspace{.2 in}

{\Large Varsha Daftardar-Gejji\\}
\textit{Department of Mathematics, Savitribai Phule Pune University,\\ Pune - 411007, India\\ vsgejji@unipune.ac.in, vsgejji@gmail.com}\\

\vskip .2in

{\bf \large Abstract}
\end{center}
\vspace{.1in}
Developing analytical methods for solving fractional partial differential equations (FPDEs) is  an active area of research. Especially finding exact solutions of FPDEs is a challenging task.  In the present paper we extend Sumudu transform iterative method (STIM) to solve a variety of time and space FPDEs as well as systems of them. We demonstrate the utility of the method by finding exact solutions to a large number of FPDEs.

\section{Introduction}
Nonlinear fractional partial differential equations (FPDEs) play important role in science and technology as they describe various nonlinear phenomena especially dealing with memory. To obtain physical information and deeper insights into the physical
aspects of the problems one has to find their exact solutions which usually
is a difficult task. For solving linear FPDEs, integral transform methods are extended successfully \cite{podlubny}. 
Various decomposition methods have been developed for solving the linear and nonlinear FPDEs such as Adomian decomposition method (ADM) \cite{adomian1994}, Homotopy perturbation method (HPM) \cite{he1999homotopy}, Daftardar-Gejji and Jafari method (DGJM) \cite{daftardar2006iterative}, the variational iteration method (VIM) \cite{he1999variational} and so on. Further, combinations of integral transforms and decomposition methods have proven to be useful. A combination of Laplace transform and DGJM (Iterative Laplace transform method (ILTM)) has been developed by Jafari {\it{et al}} \cite{jafari2013new}. A combination of HPM and Sumudu transform yields homotopy perturbation Sumudu transform method (HPSTM)\cite{singh2011homotopy}. Similarly, a  combination of Sumudu transform and ADM termed as Sumudu decomposition method (SDM) has been developed \cite{kumar2012sumudu}.
Recently, Sumudu transform iterative method (STIM) which is a combination of Sumudu transform and DGJM has been introduced and applied for solving time-fractional Cauchy reaction-diffusion equation \cite{wang2016new}.
Further, a fractional model of nonlinear Zakharov-Kuznetsov equations also have been solved using STIM \cite{prakash2018new}.
 
In this paper we extend STIM to solve time and space FPDEs as well as  systems of them. A variety of problems have been solved using STIM. In some cases, the STIM yields an exact solutions of the time and space FPDEs as well as systems of them which can be expressed in terms of the well-known Mittag-Leffler functions or fractional trigonometric functions.

The organization of this paper as follows: In section \ref{basicdefn}, we give basic definitions related to fractional calculus and Sumudu transform. In section \ref{stimeqn}, we extend STIM for time and space FPDEs. In section \ref{stimeqnexample}, we apply extended STIM to solve various time and space FPDEs. Further, in section \ref{stimsys} we extend STIM for system of time and space FPDEs. In section \ref{stimsysexample}, we apply extended STIM for system of time and space FPDEs. Conclusions are summarized in section \ref{conclusion}. 

\section{Preliminaries and Notations}\label{basicdefn}
In this section, we give some basic definitions, notations and properties of the fractional calculus (\cite{podlubny,kilbas}), which are used further in this paper.
\begin{defn}
Riemann-Liouville fractional integral of order $\alpha > 0$, of a real valued function  $f(t)$ is defined as
\begin{equation}
I_{t}^{\alpha}f(t)=\frac{1}{\Gamma(\alpha)}\int_{0}^{t}(t-s)^{\alpha-1}f(s)ds.
\end{equation}
\end{defn}
\begin{defn} 
Caputo derivative of order $\alpha> 0$ $ (n-1<\alpha<n), n\in \mathbb{N}$ of a real valued function  $f(t)$ is defined as
\begin{eqnarray}
D_t^\alpha f(t) &=& I_{t}^{n-\alpha}\Big[\frac{d^n f(t)}{dt^n}\Big], \nonumber \\ 
&=& \left\{\begin{array}{rcl}
\frac{1}{\Gamma{(n-\alpha)}}\int_0^t (t-s)^{n-\alpha-1} \frac{d^n f(s)}{ds^n} ds, & n-1<\alpha<n, \\
\frac{d^n f(t)}{dt^n}, & \alpha = n.
\end{array}\right.
\end{eqnarray}
\end{defn}
{\bf{Note:}}
\begin{enumerate}
\item  $\frac{d^{\alpha} C}{d t^{\alpha}}=0$, where $C$ is a constant.\\
\item For $\lceil \alpha \rceil = n,~n \in \mathbb{N}$, 
\begin{align}\label{rloperator}
\frac{d^{\alpha}t^p}{dt^{\alpha}}: &=
\left\{
\begin{array}{ll}
\ 0, & if ~p \in {0,1,2,...,n-1},\\
\ \frac{\Gamma(p+1)}{\Gamma(p-\alpha+1)}t^{p-\alpha}, & if ~ p \in \mathbb{N} ~and~ p \geq n,~ or~ p \neq \mathbb{N}~ and~ p > n-1.
\end{array}
\right.
\end{align}
\end{enumerate}

\begin{defn}
Riemann-Liouville time-fractional integral of order $\alpha > 0$, of a real valued function  $u(x,t)$ is defined as
\begin{equation}
I_{t}^{\alpha}u(x,t)=\frac{1}{\Gamma(\alpha)}\int_{0}^{t}(t-s)^{\alpha-1}u(x,s)ds.
\end{equation}
\end{defn}

\begin{defn}
The Caputo time-fractional derivative operator of order 
$\alpha>0$ $ (m-1<\alpha<m), n\in \mathbb{N}$ of a real valued function $u(x,t)$ is defined as 
\begin{eqnarray}
\frac{\partial^\alpha u(x,t)}{\partial t^{\alpha}} &=& I_{t}^{m-\alpha} \Big[\frac{\partial^m u(x,t)}{\partial t^m}\Big], \nonumber \\ 
&=& \left\{\begin{array}{rcl}
\frac{1}{\Gamma{(m-\alpha)}}\int_0^t (t-y)^{m-\alpha-1} \frac{\partial^m u(x,y)}  {\partial y^m} dy, & m-1<\alpha<m, \\
\frac{\partial^m u(x,t)}{\partial t^m}, & \alpha = m.
\end{array}\right.
\end{eqnarray}
\end{defn} 
Similarly, the Caputo space-fractional derivative operator $\frac{\partial^\beta u(x,t)}{\partial x^{\beta}}$ of order $\beta>0~(m-1<\beta<m), m\in \mathbb{N}$ can be defined. 

Note that: In the present paper fractional derivative $\frac{\partial^{l\beta} u(x,t)}{\partial x^{l\beta}}$, $l\in \mathbb{N}$ is taken as the sequential fractional derivative \cite{miller1993introduction} i.e.\\
\begin{eqnarray}
\frac{\partial^{l\beta}u}{\partial x^{l \beta}}=\underbrace{\frac{\partial^{\beta}}{\partial x^{\beta}}\frac{\partial^{\beta}}{\partial x^{\beta}}...\frac{\partial^{\beta}u}{\partial x^{\beta}}}_{l-times}
\end{eqnarray}

\begin{defn}
Mittag-Leffler function with two parameters $\alpha$ and $\beta$ is defined as
\begin{equation}\label{mittag2}
E_{\alpha,\beta}(z)=\sum_{k=0}^\infty \frac{z^k}{\Gamma(\alpha k+\beta)},~~~   Re(\alpha)>0, z, \beta \in \mathbb{C}.
\end{equation}
\end{defn}

{\bf{Note that:}}
\begin{enumerate}
\item The $\alpha-$th order Caputo derivative of $E_{\alpha}(a t^{\alpha})$ is  
\begin{equation}
\frac{d^{\alpha}}{d t^{\alpha}}E_{\alpha}(a t^{\alpha}) = a E_{\alpha}(a t^{\alpha}), ~\alpha > 0,~~ a \in \mathbb{R}.
\end{equation}

\item  Generalized fractional trigonometric functions for $\lceil \alpha \rceil = n$ are defined as \cite{bonilla2007fractional}
\begin{equation}
\left.
\begin{aligned}
    \cos_{\alpha}(\lambda t^{\alpha} ) &= \Re [E_{\alpha}(i\lambda^{\alpha})]=\sum_{k=0}^{\infty}\frac{(-1)^k \lambda ^{2 k} t^{(2k)\alpha}}{\Gamma(2 k \alpha+1)},   \quad\\ 
    \sin_{\alpha}(\lambda t^{\alpha} ) &= \Im [E_{\alpha}(i\lambda^{\alpha})] = \sum_{k=0}^{\infty}\frac{(-1)^k \lambda ^{2 k+1} t^{(2k+1)\alpha}}{\Gamma((2 k+1) \alpha+1)}.
\end{aligned}
\right\}
\end{equation}

\item The Caputo derivative of fractional trigonometric functions are defined as
\begin{equation}
\left.
\begin{aligned}
\frac{d^{\alpha}}{d t^{\alpha}}\cos_{\alpha}(\lambda t^{\alpha})  &=  -\lambda \sin_{\alpha}(\lambda t^{\alpha}), \quad\\ 
  \frac{d^{\alpha}}{d t^{\alpha}}\sin_{\alpha}(\lambda t^{\alpha}) &= \lambda \cos_{\alpha}(\lambda t^{\alpha}).
\end{aligned}
\right\}
\end{equation}
\end{enumerate}

\begin{defn}
\cite{belgacem2006sumudu} The Sumudu transform over the set of functions \\$ A = \{f(t)~ | ~\exists~ M, \tau_1, \tau_2 >0,$ such that$~ |f(x,t)|< M e^{|t|/\tau_j} ~if~\\  t \in (-1)^j \times [0,\infty)\}$ is defined as 
\begin{eqnarray}
S[f(t)] =F(\omega)= \int_{0}^{\infty} e^{-t} f(\omega t)dt, ~~ \omega \in (-\tau_1, \tau_2).
\end{eqnarray}
\end{defn}
One of the basic property of Sumudu transform is 
\begin{eqnarray}
S\Big[\frac{t^{\alpha}}{\Gamma(\alpha+1)}\Big]= \omega^{\alpha}, ~~ \alpha > 0.
\end{eqnarray}
Sumudu inverse transform of $\omega^{\alpha}$ is defined as
\begin{eqnarray}
S^{-1}[ \omega^{\alpha}]=\frac{t^{\alpha}}{\Gamma(\alpha+1)}, ~~ \alpha > 0.
\end{eqnarray}

\begin{defn}
\cite{amer2018solving} The Sumudu transform of Caputo time-fractional derivative of $f(x,t)$ of order $\gamma>0$ is defined as 
\begin{equation}\label{sproperty}
\left.
\begin{aligned}
 S\Big[\frac{\partial^{\gamma}f(x,t)}{\partial t^{\gamma}}\Big] &= \omega^{-\gamma} S[f(x,t)]-\sum_{k=0}^{m-1}\Big[\omega^{-\gamma+k}\frac{\partial^k f(x,0)}{\partial t^k}\Big], \quad\\ 
     &~~~~~~~~~~~~~~~~~~~~~~~~  m-1 < \gamma \leq m, ~ m \in \mathbb{N}.
\end{aligned}
\right\}
\end{equation}
\end{defn}

\section{STIM for time and space FPDEs}\label{stimeqn}
In this section, we extend STIM \cite{wang2016new} for solving time and space FPDEs.\\We consider the following general time and space FPDE: 
\begin{equation}\label{steq}
\left.
\begin{aligned}
 \frac{\partial^{\gamma}u}{\partial t^{\gamma}} &= \mathscr{F}\Big(x, u,\frac{\partial^{\beta}u}{\partial x^{\beta}},...,\frac{\partial^{l \beta}u}{\partial x^{l \beta} } \Big),~  m-1<\gamma \leq m, \quad\\ 
     &~~~~~~~~~~~~~~~~~~~ n-1<\beta \leq n,~ l, m, n \in \mathbb{N},
\end{aligned}
\right\}
\end{equation}

along with the initial conditions
\begin{eqnarray}\label{ini}
\frac{\partial^k u(x,0)}{\partial t^k}=h_ k(x), ~~k=0,1,2,...,m-1,
\end{eqnarray}
where $\mathscr{F}\Big(x, u,\frac{\partial^{\beta}u}{\partial x^{\beta}},...,\frac{\partial^{l \beta}u}{\partial x^{l \beta} } \Big) $ is a linear/nonlinear operator and $u = u(x,t)$ is the unknown function.

Taking the Sumudu transform of both sides of Eq. \eqref{steq} and simplifying, we get
\begin{eqnarray}\label{singles}
S[u(x,t)]=\sum_{k=0}^{m-1} \Big[\omega^{k}\frac{\partial^k u(x,0)}{\partial t^k}\Big]+\omega^{\gamma}S\Big[\mathscr{F}\Big(x, u,\frac{\partial^{\beta}u}{\partial x^{\beta}},...,\frac{\partial^{l \beta}u}{\partial x^{l \beta} } \Big)\Big].
\end{eqnarray}
The inverse Sumudu transform of Eq. \eqref{singles} leads to
\begin{eqnarray}\label{stn}
u(x,t) =S^{-1}\Big(\sum_{k=0}^{m-1} \Big[\omega^{k}\frac{\partial^k u(x,0)}{\partial t^k}\Big]\Big) + S^{-1} \Big[\omega^{\gamma}S\Big(\mathscr{F}\Big(x, u,\frac{\partial^{\beta}u}{\partial x^{\beta}},...,\frac{\partial^{l \beta}u}{\partial x^{l \beta} } \Big)\Big)\Big].
\end{eqnarray}
Eq. \eqref{stn} can be written as 
\begin{eqnarray}\label{ufn}
u(x,t) = f(x,t)+N\Big(x, u,\frac{\partial^{\beta}u}{\partial x^{\beta}},...,\frac{\partial^{l \beta}u}{\partial x^{l \beta} } \Big),
\end{eqnarray}
where 
\begin{equation}\label{asume}
\left.
\begin{aligned}
    f(x,t) &= S^{-1}\Big(\sum_{k=0}^{m-1}\Big[ \omega^{k}\frac{\partial^k u(x,0)}{\partial t^k}\Big]\Big), \quad\\ 
    N\Big(x, u,\frac{\partial^{\beta}u}{\partial x^{\beta}},...,\frac{\partial^{l \beta}u}{\partial x^{l \beta} } \Big) &= S^{-1} \Big[\omega^{\gamma}S\Big(\mathscr{F}\Big(x, u,\frac{\partial^{\beta}u}{\partial x^{\beta}},...,\frac{\partial^{l \beta}u}{\partial x^{l \beta} } \Big)\Big)\Big],
\end{aligned}
\right\}
\end{equation}
here $f$ is known function and $N$ is a linear/nonlinear operator. 

Functional equations of the form \eqref{ufn} can be solved by the DGJ decomposition method introduced by Daftardar-Gejji and Jafari \cite{daftardar2006iterative}.

DGJM represents the solution as an infinite series:
\begin{equation}\label{ssum}
u = \sum_{i=0}^{\infty} u_i,
\end{equation}
where the terms $u_i$ are calculated recursively. 
The operator $N$ can be decomposed as 
\begin{eqnarray}
N\Big(x, \sum_{i=0}^{\infty}u_i,\frac{\partial^{\beta}(\sum_{i=0}^{\infty}u_i)}{\partial x^{\beta}},...,\frac{\partial^{l \beta}(\sum_{i=0}^{\infty}u_i)}{\partial x^{l \beta} } \Big) =  N\Big(x, u_0,\frac{\partial^{\beta}u_0}{\partial x^{\beta}},...,\frac{\partial^{l \beta}u_0}{\partial x^{l \beta} } \Big) \nonumber \\
+\sum_{j=1}^\infty \Big(N\Big(x, \sum_{i=0}^{j}u_i,\frac{\partial^{\beta}(\sum_{i=0}^{j}u_i)}{\partial x^{\beta}},...,\frac{\partial^{l \beta}(\sum_{i=0}^{j}u_i)}{\partial x^{l \beta} } \Big)\Big) \nonumber \\ -\sum_{j=1}^\infty \Big(N\Big(x, \sum_{i=0}^{j-1}u_i,\frac{\partial^{\beta}(\sum_{i=0}^{j-1}u_i)}{\partial x^{\beta}},...,\frac{\partial^{l \beta}(\sum_{i=0}^{j-1}u_i)}{\partial x^{l \beta} } \Big) \Big).
\end{eqnarray}
\begin{equation}\label{decom}
\left.
\begin{aligned}
   S^{-1} \Big[\omega^{\gamma}S\Big(\mathscr{F}\Big(x,\sum_{i=0}^{\infty}u_i ,\frac{\partial^{\beta}(\sum_{i=0}^{\infty} u_i)}{\partial x^{\beta}}, ...,\frac{\partial^{l \beta}(\sum_{i=0}^{\infty} u_i)}{\partial x^{l \beta} } \Big)\Big)\Big] & \quad\\ 
     = S^{-1} \Big[\omega^{\gamma}S\Big(\mathscr{F}\Big(x, u_0,\frac{\partial^{\beta}u_0}{\partial x^{\beta}}, ...,\frac{\partial^{l \beta}u_0}{\partial x^{l \beta} }\Big)\Big)\Big]& \quad \\
    +\sum_{j=1}^{\infty} S^{-1} \Big[\omega^{\gamma}S\Big(\mathscr{F}\Big(x,\sum_{i=0}^{j}u_i ,\frac{\partial^{\beta}(\sum_{i=0}^{j} u_i)}{\partial x^{\beta}}, ...,\frac{\partial^{l \beta}(\sum_{i=0}^{j} u_i)}{\partial x^{l \beta} }\Big)\Big)\Big]&\quad \\
     -\sum_{j=1}^{\infty} S^{-1} \Big[\omega^{\gamma}S\Big(\mathscr{F}\Big(x, \sum_{i=0}^{j-1}u_i ,\frac{\partial^{\beta}(\sum_{i=0}^{j-1} u_i)}{\partial x^{\beta}},...,\frac{\partial^{l \beta}(\sum_{i=0}^{j-1} u_i)}{\partial x^{l \beta} }\Big)\Big)\Big].
\end{aligned}
\right\}
\end{equation}

Using Eqns. (\ref{ssum},\ref{decom}) in Eq. \eqref{ufn}, we get
\begin{eqnarray}
\sum_{i=0}^{\infty} u_i = S^{-1}\Big(\sum_{k=0}^{m-1} \Big[\omega^{k}\frac{\partial^k u(x,0)}{\partial t^k}\Big]\Big)+S^{-1} \Big[\omega^{\gamma}S\Big(\mathscr{F}\Big(x, u_0,\frac{\partial^{\beta}u_0}{\partial x^{\beta}}, ...,\frac{\partial^{l \beta}u_0}{\partial x^{l \beta} }\Big)\Big)\Big] \nonumber \\
+\sum_{j=1}^{\infty} \Big(S^{-1} \Big[\omega^{\gamma}S\Big(\mathscr{F}\Big(x, \sum_{i=0}^{j}u_i ,\frac{\partial^{\beta}(\sum_{i=0}^{j} u_i)}{\partial x^{\beta}},...,\frac{\partial^{l \beta}(\sum_{i=0}^{j} u_i)}{\partial x^{l \beta} }\Big)\Big)\Big] \nonumber \\
- S^{-1} \Big[\omega^{\gamma}S\Big(\mathscr{F}\Big(x, \sum_{i=0}^{j-1}u_i ,\frac{\partial^{\beta}(\sum_{i=0}^{j-1} u_i)}{\partial x^{\beta}}, ...,\frac{\partial^{l \beta}(\sum_{i=0}^{j-1} u_i)}{\partial x^{l \beta} }\Big)\Big)\Big]\Big).
\end{eqnarray}
We define the recurrence relation as follows:

\begin{equation}\label{nstim}
\left.
\begin{aligned}
u_0 &= S^{-1}\Big(\sum_{k=0}^{m-1} \Big[\omega^{k}\frac{\partial^k u(x,0)}{\partial t^k}\Big]\Big),\\ 
u_1 &= S^{-1} \Big[\omega^{\gamma}S\Big(\mathscr{F}\Big(x, u_0,\frac{\partial^{\beta}u_0}{\partial x^{\beta}}, ...,\frac{\partial^{l \beta}u_0}{\partial x^{l \beta} }\Big)\Big)\Big],\\
\hspace*{-2cm}u_{r+1}&= 
S^{-1} \Big[\omega^{\gamma}S\Big(\mathscr{F}\Big(x, \sum_{i=0}^{r}u_i ,\frac{\partial^{\beta}(\sum_{i=0}^{r} u_i)}{\partial x^{\beta}}, ...,\frac{\partial^{l \beta}(\sum_{i=0}^{r} u_i)}{\partial x^{l \beta} }\Big)\Big)\Big]\\
&- S^{-1} \Big[\omega^{\gamma}S\Big(\mathscr{F}\Big(x, \sum_{i=0}^{r-1}u_i ,\frac{\partial^{\beta}(\sum_{i=0}^{r-1} u_i)}{\partial x^{\beta}}, ...,\frac{\partial^{l \beta}(\sum_{i=0}^{r-1} u_i)}{\partial x^{l \beta} }\Big)\Big)\Big],\\&~~~~~~~~~~for ~~~~ r\geq 1.
\end{aligned}
\right\}
\end{equation}

The r-term approximate solution of Eqns. (\ref{steq}-\ref{ini}) is given by\\ $u \approx u_0+u_1+ \cdots +u_{r-1}$. For the convergence of DGJM we refer the reader to \cite{bhalekar2011convergence}.

\section{Illustrative Examples}\label{stimeqnexample}
In this section, we solve various linear and nonlinear time and space FPDEs using STIM derived in \ref{stimeqn}.

\begin{ex}
Consider the following time and space linear fractional Newell-Whitehead- Segel equation
\end{ex}
\begin{eqnarray}\label{nw}
\frac{\partial^{\alpha}u}{\partial t^{\alpha}}=\frac{\partial^{2 \beta}u}{\partial x^{2 \beta}}-3 u,~~ t>0, ~\alpha, \beta \in (0,1],
\end{eqnarray}
along with the initial condition
\begin{eqnarray}\label{nwi}
u(x,0)=E_{\beta}(2 x^{\beta}).
\end{eqnarray}

Taking the Sumudu transform of both sides of  Eq. \eqref{nw}
\begin{eqnarray}
S\Big[\frac{\partial^{\alpha}u}{\partial t^{\alpha}} \Big]=S\Big[\frac{\partial^{2 \beta}u}{\partial x^{2 \beta}}-3 u \Big].
\end{eqnarray}
Using the property of Sumudu transform \eqref{sproperty} 
\begin{eqnarray}\label{nws}
S[u(x,t)]=u(x,0)+\omega^{\alpha}S\Big[\frac{\partial^{2 \beta}u}{\partial x^{2 \beta}}-3 u \Big].
\end{eqnarray}
Taking the inverse Sumudu transform of both sides of Eq. \eqref{nws}
\begin{eqnarray}
u(x,t)=S^{-1}[u(x,0)]+S^{-1}\Big(\omega^{\alpha}S\Big[\frac{\partial^{2 \beta}u}{\partial x^{2 \beta}}-3 u \Big] \Big).
\end{eqnarray}

In view of the recurrence relation \eqref{nstim}, we get 
\begin{eqnarray}
u_0 &=& S^{-1}[u(x,0)]=E_{\beta}(2 x^{\beta}),\nonumber \\
u_1 &=& S^{-1}\Big(\omega^{\alpha}S\Big[\frac{\partial^{2 \beta}u_0}{\partial x^{2 \beta}}-3 u_0 \Big] \Big)=\frac{E_{\beta}(2 x^{\beta}) t^{\alpha }}{\Gamma (\alpha +1)}, \nonumber \\
u_2 &=& \frac{E_{\beta}(2 x^{\beta}) t^{2\alpha }}{\Gamma (2\alpha +1)}, \nonumber \\
u_3 &=& \frac{E_{\beta}(2 x^{\beta}) t^{3\alpha }}{\Gamma (3\alpha +1)},\nonumber \\
\vdots
\end{eqnarray}

Hence, the series form solution of Eqns. (\ref{nw}-\ref{nwi}) is given by
\begin{eqnarray}
u(x,t) &=& u_0+u_1+u_2+...\nonumber \\
&=& E_{\beta}(2 x^{\beta})+\frac{E_{\beta}(2 x^{\beta}) t^{\alpha }}{\Gamma (\alpha +1)}+\frac{E_{\beta}(2 x^{\beta}) t^{2\alpha }}{\Gamma (2\alpha +1)}+\frac{E_{\beta}(2 x^{\beta}) t^{3\alpha }}{\Gamma (3\alpha +1)}+....,\nonumber \\
\end{eqnarray}
which converges to 
\begin{eqnarray}
u(x,t)&=& E_{\beta}(2 x^{\beta})E_{\alpha}(t^{\alpha}).
\end{eqnarray}

\begin{ex}
Consider the following time and space linear fractional diffusion equation:
\begin{eqnarray}\label{diff}
\frac{\partial^{\alpha}u}{\partial t^{\alpha}} = c \frac{\partial^{\beta}u}{\partial x^{\beta}},~~t>0, ~\alpha, \beta \in (0,1],
\end{eqnarray}
along with the initial condition
\begin{eqnarray}\label{diffi}
u(x,0)=a+bx^{\beta}, ~~a,b \in \mathbb{R}.
\end{eqnarray}
\end{ex}
Taking the Sumudu transform of both sides of Eq. \eqref{diff}, we get
\begin{eqnarray}
S\Big[\frac{\partial^{\alpha}u}{\partial t^{\alpha}}\Big] = c S\Big[ \frac{\partial^{\beta}u}{\partial x^{\beta}}\Big].
\end{eqnarray}
Using the property of Sumudu transform \eqref{sproperty}, we get
\begin{eqnarray}\label{diffs}
S[u(x,t)] = u(x,0) + c \omega^{\alpha}\Big[\frac{\partial^{\beta}u}{\partial x^{\beta}} \Big]. 
\end{eqnarray}
Taking the inverse Sumudu transform  of both sides of Eq. \eqref{diffs}
\begin{eqnarray}
u(x,t) = S^{-1}[u(x,0)] + c S^{-1}\Big( \omega^{\alpha}\Big[\frac{\partial^{\beta}u}{\partial x^{\beta}} \Big]\Big). 
\end{eqnarray}
Using the recurrence relation \eqref{nstim}, we get
\begin{eqnarray}
u_0 &=& S^{-1}[u(x,0)] = a+bx^{\beta}, \nonumber \\ 
u_1 &=& c S^{-1}\Big( \omega^{\alpha}\Big[\frac{\partial^{\beta}u_0}{\partial x^{\beta}} \Big]\Big), \nonumber \\
&=& c b\Gamma(\beta+1) \frac{t^{\alpha}}{\Gamma(\alpha+1)}, \nonumber \\
u_i &=& 0,~~ i\geq 2.
\end{eqnarray}

Hence, we obtain the exact solution of Eqns. (\ref{diff}-\ref{diffi} as
\begin{eqnarray}
u(x,t)=\sum_{i=0}^{\infty}u_i = a+\frac{b c \Gamma (\beta +1) t^{\alpha }}{\Gamma (\alpha +1)}+b x^{\beta }.
\end{eqnarray}

\begin{ex}
Consider the following time and  space fractional equation
\begin{eqnarray}\label{frac2}
\frac{\partial^{\alpha}u}{\partial t^{\alpha}} &=& \Big( \frac{\partial^{\beta}u}{\partial x^{\beta}} \Big)^2 - u \Big( \frac{\partial^{\beta}u}{\partial x^{\beta}} \Big), ~~t>0,~ \alpha, \beta \in (0,1],
\end{eqnarray}
along with the initial condition 
\begin{eqnarray}\label{frac2i}
u(x,0) &=& 3+ \frac{5}{2} E_{\beta}(x^{\beta}).
\end{eqnarray}
\end{ex}
Taking the Sumudu transform of both sides of Eq. \eqref{frac2}, we get
\begin{eqnarray}
S\Big[\frac{\partial^{\alpha}u}{\partial t^{\alpha}}\Big] = S\Big[\Big( \frac{\partial^{\beta}u}{\partial x^{\beta}} \Big)^2 - u \Big( \frac{\partial^{\beta}u}{\partial x^{\beta}} \Big)\Big].
\end{eqnarray}
Using the property of Sumudu transform \eqref{sproperty}, we get
\begin{eqnarray}\label{frac2s}
S[u(x,t)] = u(x,0)+\omega^{\alpha} \Big( S\Big[\Big( \frac{\partial^{\beta}u}{\partial x^{\beta}} \Big)^2 - u \Big( \frac{\partial^{\beta}u}{\partial x^{\beta}} \Big)\Big]\Big).
\end{eqnarray}
Now taking the inverse Sumudu transform of both sides of Eq. \eqref{frac2s}
\begin{eqnarray}
u(x,t) = S^{-1}[u(x,0)]+S^{-1} \Big(\omega^{\alpha} \Big( S\Big[\Big( \frac{\partial^{\beta}u}{\partial x^{\beta}} \Big)^2 - u \Big( \frac{\partial^{\beta}u}{\partial x^{\beta}} \Big)\Big]\Big)\Big).
\end{eqnarray}

Using the recurrence relation \eqref{nstim}
\begin{eqnarray}
u_0 &=& S^{-1}[u(x,0)] = 3+ \frac{5}{2} E_{\beta}(x^{\beta}), \nonumber \\
u_1 &=& S^{-1} \Big(\omega^{\alpha} \Big( S\Big[\Big( \frac{\partial^{\beta}u_0}{\partial x^{\beta}} \Big)^2 - u_0 \Big( \frac{\partial^{\beta}u_0}{\partial x^{\beta}} \Big)\Big]\Big)\Big) = -\frac{15 t^{\alpha } E_{\beta }(x^{\beta})}{2 \Gamma (\alpha +1)}, \nonumber \\
u_2 &=& S^{-1} \Big(\omega^{\alpha} \Big( S\Big[\Big( \frac{\partial^{\beta}(u_0+u_1)}{\partial x^{\beta}} \Big)^2 - (u_0+u_1) \Big( \frac{\partial^{\beta}(u_0+u_1)}{\partial x^{\beta}} \Big)\Big]\Big)\Big), \nonumber \\
&-& S^{-1} \Big(\omega^{\alpha} \Big( S\Big[\Big( \frac{\partial^{\beta}u_0}{\partial x^{\beta}} \Big)^2 - u_0 \Big( \frac{\partial^{\beta}u_0}{\partial x^{\beta}} \Big)\Big]\Big)\Big), \nonumber \\
&=& \frac{45 t^{2 \alpha } E_{\beta }(x^{\beta})}{2 \Gamma (2 \alpha +1)}, \nonumber \\
u_3 &=& - \frac{135 t^{3 \alpha } E_{\beta }(x^{\beta})}{2 \Gamma (3 \alpha +1)}, \nonumber \\
u_4 &=& \frac{405 t^{4 \alpha } E_{\beta }(x^{\beta})}{2 \Gamma (4 \alpha +1)}, \nonumber \\
\vdots
\end{eqnarray}
Hence, the series solution of Eq. $(\ref{frac2})$ along with the initial condition (\ref{frac2i}) is given by 
\begin{eqnarray}
u(x,t)&=& 3+ \frac{5}{2} E_{\beta}(x^{\beta})-\frac{15 t^{\alpha } E_{\beta }(x^{\beta})}{2 \Gamma (\alpha +1)}+\frac{45 t^{2 \alpha } E_{\beta }(x^{\beta})}{2 \Gamma (2 \alpha +1)}- \frac{135 t^{3 \alpha } E_{\beta }(x^{\beta})}{2 \Gamma (3 \alpha +1)} \nonumber \\
&+&\frac{405 t^{4 \alpha } E_{\beta }(x^{\beta})}{2 \Gamma (4 \alpha +1)}-\cdots.
\end{eqnarray}
This leads to the following  closed form solution:
\begin{eqnarray}
u(x,t) = 3+\Big[\frac{5}{2}E_{\alpha}(-3t^{\alpha})\Big]E_{\beta}(x^{\beta}),
\end{eqnarray} 
which is the same as obtained in \cite{choudhary2017invariant}.

\begin{ex}
Consider the following time space fractional heat equation:
\begin{eqnarray}\label{heat}
\frac{\partial^{\alpha}u}{\partial t^{\alpha}}= \frac{\partial^{\beta}}{\partial x^{\beta}}\Big(u \frac{\partial^{\beta}u}{\partial x^{\beta}} \Big),  ~~t>0,~ \alpha, \beta \in (0,1],
\end{eqnarray}
with the initial condition 
\begin{eqnarray}\label{heati}
u(x,0)= a+b x^{\beta}, ~~ a, b \in \mathbb{R}.
\end{eqnarray}
\end{ex}
Taking the Sumudu transform of both sides of Eq. \eqref{heat}, we get
\begin{eqnarray}\label{heats}
S\Big[\frac{\partial^{\alpha}u}{\partial t^{\alpha}}\Big]=S\Big[ \frac{\partial^{\beta}}{\partial x^{\beta}}\Big(u \frac{\partial^{\beta}u}{\partial x^{\beta}} \Big)\Big], \nonumber \\
\implies S[u(x,t)]= u(x,0)+ \omega^{\alpha}\Big(S\Big[ \frac{\partial^{\beta}}{\partial x^{\beta}}\Big(u \frac{\partial^{\beta}u}{\partial x^{\beta}} \Big)\Big] \Big), \nonumber \\
\implies u(x,t)=S^{-1}[u(x,0)]+S^{-1} \Big(\omega^{\alpha}\Big(S\Big[ \frac{\partial^{\beta}}{\partial x^{\beta}}\Big(u \frac{\partial^{\beta}u}{\partial x^{\beta}} \Big)\Big] \Big) \Big).
\end{eqnarray}
Now using the recurrence relation \eqref{nstim}
\begin{eqnarray}
u_0 &=& S^{-1}[u(x,0)] = a+bx^{\beta}, \nonumber \\
u_1 &=& S^{-1} \Big(\omega^{\alpha}\Big(S\Big[ \frac{\partial^{\beta}}{\partial x^{\beta}}\Big(u_0 \frac{\partial^{\beta}u_0}{\partial x^{\beta}} \Big)\Big] \Big) \Big),\nonumber \\
 &=& b^2 (\Gamma(\beta+1))^2 \frac{t^{\alpha}}{\Gamma(\alpha+1)}, \nonumber \\
u_i &=& 0,~~ \forall~ i \geq 2.
\end{eqnarray}
Hence, the solution turns out to be:
\begin{eqnarray}
u(x,t) &=& a+bx^{\beta}+b^2 (\Gamma(\beta+1))^2 \frac{t^{\alpha}}{\Gamma(\alpha+1)}.
\end{eqnarray}

\begin{ex}
Consider the following time and space fractional thin film equation
\begin{eqnarray}\label{thin}
\frac{\partial^{\alpha}u}{\partial t^{\alpha}}=-u\Big(\frac{\partial^{4 \beta}u}{\partial x^{4\beta}} \Big) + \eta \Big(\frac{\partial^{\beta}u}{\partial x^{\beta}}\Big)\Big(\frac{\partial^{3 \beta}u}{\partial x^{3 \beta}} \Big)+ \zeta \Big(\frac{\partial^{2 \beta}u}{\partial x^{2\beta}} \Big)^2, ~ t>0~ \alpha, \beta \in (0,1],
\end{eqnarray}
along with the initial condition
\begin{eqnarray}\label{thini}
u(x,0)=a+b x^{\beta }+c x^{2 \beta }+d x^{3 \beta }, ~~a,b,c,d \in \mathbb{R}.
\end{eqnarray}
\end{ex}
Taking the Sumudu transform of both sides of Eq. \eqref{thin}, we get
\begin{eqnarray}
S\Big[\frac{\partial^{\alpha}u}{\partial t^{\alpha}}\Big]=S\Big[-u\Big(\frac{\partial^{4 \beta}u}{\partial x^{4\beta}} \Big) + \eta \Big(\frac{\partial^{\beta}u}{\partial x^{\beta}}\Big)\Big(\frac{\partial^{3 \beta}u}{\partial x^{3 \beta}} \Big)+ \zeta \Big(\frac{\partial^{2 \beta}u}{\partial x^{2\beta}} \Big)^2\Big].
\end{eqnarray}
After simplification, we get
\begin{eqnarray}
\hspace*{-0.5cm}u(x,t)&=& S^{-1}[u(x,0)] \nonumber \\ 
&+& S^{-1}\Big(\omega^{\alpha} S\Big[-u\Big(\frac{\partial^{4 \beta}u}{\partial x^{4\beta}} \Big) + \eta \Big(\frac{\partial^{\beta}u}{\partial x^{\beta}}\Big)\Big(\frac{\partial^{3 \beta}u}{\partial x^{3 \beta}} \Big)+ \zeta \Big(\frac{\partial^{2 \beta}u}{\partial x^{2\beta}} \Big)^2\Big]\Big).
\end{eqnarray}
In view of the recurrence relation \eqref{nstim},
\begin{eqnarray}
\hspace*{-0.5cm}u_0 &=& S^{-1}[u(x,0)] = a+b x^{\beta }+c x^{2 \beta }+d x^{3 \beta }, \nonumber\\
u_1 &=& S^{-1}\Big(\omega^{\alpha} S\Big[-u_0\Big(\frac{\partial^{4 \beta}u_0}{\partial x^{4\beta}} \Big) + \eta \Big(\frac{\partial^{\beta}u_0}{\partial x^{\beta}}\Big)\Big(\frac{\partial^{3 \beta}u_0}{\partial x^{3 \beta}} \Big)+ \zeta \Big(\frac{\partial^{2 \beta}u_0}{\partial x^{2\beta}} \Big)^2\Big]\Big), \nonumber \\
&=& \frac{b \eta d \Gamma (\beta +1) \Gamma (3 \beta +1) t^{\alpha }}{\Gamma (\alpha +1)}+\frac{\eta c d \Gamma (2 \beta +1) \Gamma (3 \beta +1) t^{\alpha } x^{\beta }}{\Gamma (\alpha +1) \Gamma (\beta +1)}\nonumber \\
&+&\frac{\eta d^2 \Gamma (3 \beta +1)^2 t^{\alpha } x^{2 \beta }}{\Gamma (\alpha +1) \Gamma (2 \beta +1)}+\frac{c^2 \zeta  \Gamma (2 \beta +1)^2 t^{\alpha }}{\Gamma (\alpha +1)}\nonumber \\ 
&+&\frac{2 c \zeta  d \Gamma (2 \beta +1) \Gamma (3 \beta +1) t^{\alpha } x^{\beta }}{\Gamma (\alpha +1) \Gamma (\beta +1)} + \frac{\zeta  d^2 \Gamma (3 \beta +1)^2 t^{\alpha } x^{2 \beta }}{\Gamma (\alpha +1) \Gamma (\beta +1)^2}. \nonumber \\
u_2 &=& \frac{\eta^2 c d^2 \Gamma (2 \beta +1) \Gamma (3 \beta +1)^2 t^{2 \alpha }}{\Gamma (2 \alpha +1)}+\frac{2 \zeta ^2 d^3 \Gamma (2 \beta +1) \Gamma (3 \beta +1)^3 t^{2 \alpha } x^{\beta }}{\Gamma (2 \alpha +1) \Gamma (\beta +1)^3}   \nonumber \\ &+&
\frac{\eta^2 d^3 \Gamma (3 \beta +1)^3 t^{2 \alpha } x^{\beta }}{\Gamma (2 \alpha +1) \Gamma (\beta +1)}+\frac{4 \eta c \zeta  d^2 \Gamma (2 \beta +1) \Gamma (3 \beta +1)^2 t^{2 \alpha }}{\Gamma (2 \alpha +1)}\nonumber \\
&+& \frac{\eta \zeta  d^3 \Gamma (2 \beta +1) \Gamma (3 \beta +1)^3 t^{2 \alpha } x^{\beta }}{\Gamma (2 \alpha +1) \Gamma (\beta +1)^3}+\frac{2 \eta \zeta  d^3 \Gamma (3 \beta +1)^3 t^{2 \alpha } x^{\beta }}{\Gamma (2 \alpha +1) \Gamma (\beta +1)} \nonumber \\ 
&+&\frac{2 c \zeta ^2 d^2 \Gamma (2 \beta +1)^2 \Gamma (3 \beta +1)^2 t^{2 \alpha }}{\Gamma (2 \alpha +1) \Gamma (\beta +1)^2}\nonumber \\ 
&+&\frac{2 \eta \zeta ^2 d^4 \Gamma (2 \alpha +1) \Gamma (2 \beta +1) \Gamma (3 \beta +1)^4 t^{3 \alpha }}{\Gamma (\alpha +1)^2 \Gamma (3 \alpha +1) \Gamma (\beta +1)^2}     
 \nonumber \\ 
&+&\frac{\eta^2 \zeta  d^4 \Gamma (2 \alpha +1) \Gamma (3 \beta +1)^4 t^{3 \alpha }}{\Gamma (\alpha +1)^2 \Gamma (3 \alpha +1)} \nonumber \\ &+&\frac{\zeta ^3 d^4 \Gamma (2 \alpha +1) \Gamma (2 \beta +1)^2 \Gamma (3 \beta +1)^4 t^{3 \alpha }}{\Gamma (\alpha +1)^2 \Gamma (3 \alpha +1) \Gamma (\beta +1)^4},\nonumber\\
u_3 &=& \frac{\eta^3 d^4 \Gamma (3 \beta +1)^4 t^{3 \alpha }}{\Gamma (3 \alpha +1)}+\frac{\eta^2 \zeta  d^4 \Gamma (2 \beta +1) \Gamma (3 \beta +1)^4 t^{3 \alpha }}{\Gamma (3 \alpha +1) \Gamma (\beta +1)^2} \nonumber \\
&+&\frac{2 \eta^2 \zeta  d^4 \Gamma (3 \beta +1)^4 t^{3 \alpha }}{\Gamma (3 \alpha +1)}+\frac{2 \eta \zeta ^2 d^4 \Gamma (2 \beta +1) \Gamma (3 \beta +1)^4 t^{3 \alpha }}{\Gamma (3 \alpha +1) \Gamma (\beta +1)^2},\nonumber \\
u_i &=& 0~~ \forall ~~i\geq 4.
\end{eqnarray}
Hence, we obtain the exact solution of Eqns. (\ref{thin}-\ref{thini}) as
\begin{eqnarray}
u(x,t) &=& a+b x^{\beta }+c x^{2 \beta }+d x^{3 \beta }\nonumber \\ 
&+&\frac{b \eta d \Gamma (\beta +1) \Gamma (3 \beta +1) t^{\alpha }}{\Gamma (\alpha +1)}+\frac{\eta c d \Gamma (2 \beta +1) \Gamma (3 \beta +1) t^{\alpha } x^{\beta }}{\Gamma (\alpha +1) \Gamma (\beta +1)} \nonumber \\ 
&+&\frac{\eta d^2 \Gamma (3 \beta +1)^2 t^{\alpha } x^{2 \beta }}{\Gamma (\alpha +1) \Gamma (2 \beta +1)}+\frac{c^2 \zeta  \Gamma (2 \beta +1)^2 t^{\alpha }}{\Gamma (\alpha +1)}\nonumber \\ 
&+& \frac{2 c \zeta  d \Gamma (2 \beta +1) \Gamma (3 \beta +1) t^{\alpha } x^{\beta }}{\Gamma (\alpha +1) \Gamma (\beta +1)} + \frac{\zeta  d^2 \Gamma (3 \beta +1)^2 t^{\alpha } x^{2 \beta }}{\Gamma (\alpha +1) \Gamma (\beta +1)^2}\nonumber \\ 
&+& \frac{\eta^2 c d^2 \Gamma (2 \beta +1) \Gamma (3 \beta +1)^2 t^{2 \alpha }}{\Gamma (2 \alpha +1)}+\frac{2 \zeta ^2 d^3 \Gamma (2 \beta +1) \Gamma (3 \beta +1)^3 t^{2 \alpha } x^{\beta }}{\Gamma (2 \alpha +1) \Gamma (\beta +1)^3} \nonumber \\ 
&+&\frac{\eta^2 d^3 \Gamma (3 \beta +1)^3 t^{2 \alpha } x^{\beta }}{\Gamma (2 \alpha +1) \Gamma (\beta +1)}+\frac{4 \eta c \zeta  d^2 \Gamma (2 \beta +1) \Gamma (3 \beta +1)^2 t^{2 \alpha }}{\Gamma (2 \alpha +1)}\nonumber \\
&+&\frac{\eta \zeta  d^3 \Gamma (2 \beta +1) \Gamma (3 \beta +1)^3 t^{2 \alpha } x^{\beta }}{\Gamma (2 \alpha +1) \Gamma (\beta +1)^3}+\frac{2 \eta \zeta  d^3 \Gamma (3 \beta +1)^3 t^{2 \alpha } x^{\beta }}{\Gamma (2 \alpha +1) \Gamma (\beta +1)}\nonumber \\
&+&\frac{2 c \zeta ^2 d^2 \Gamma (2 \beta +1)^2 \Gamma (3 \beta +1)^2 t^{2 \alpha }}{\Gamma (2 \alpha +1) \Gamma (\beta +1)^2}\nonumber \\ 
&+&\frac{2 \eta \zeta ^2 d^4 \Gamma (2 \alpha +1) \Gamma (2 \beta +1) \Gamma (3 \beta +1)^4 t^{3 \alpha }}{\Gamma (\alpha +1)^2 \Gamma (3 \alpha +1) \Gamma (\beta +1)^2}\nonumber \\ 
&+& \frac{\eta^2 \zeta  d^4 \Gamma (2 \alpha +1) \Gamma (3 \beta +1)^4 t^{3 \alpha }}{\Gamma (\alpha +1)^2 \Gamma (3 \alpha +1)}+\frac{\zeta ^3 d^4 \Gamma (2 \alpha +1) \Gamma (2 \beta +1)^2 \Gamma (3 \beta +1)^4 t^{3 \alpha }}{\Gamma (\alpha +1)^2 \Gamma (3 \alpha +1) \Gamma (\beta +1)^4} \nonumber \\ 
&+&\frac{\eta^3 d^4 \Gamma (3 \beta +1)^4 t^{3 \alpha }}{\Gamma (3 \alpha +1)}+\frac{\eta^2 \zeta  d^4 \Gamma (2 \beta +1) \Gamma (3 \beta +1)^4 t^{3 \alpha }}{\Gamma (3 \alpha +1) \Gamma (\beta +1)^2} \nonumber \\
&+&\frac{2 \eta^2 \zeta  d^4 \Gamma (3 \beta +1)^4 t^{3 \alpha }}{\Gamma (3 \alpha +1)}+\frac{2 \eta \zeta ^2 d^4 \Gamma (2 \beta +1) \Gamma (3 \beta +1)^4 t^{3 \alpha }}{\Gamma (3 \alpha +1) \Gamma (\beta +1)^2}.
\end{eqnarray}

\begin{ex}
Consider the following time and space fractional dispersive Boussinesq equation 
\begin{eqnarray}\label{db}
\frac{\partial^{2\alpha}u}{\partial t^{2\alpha}}= \frac{\partial^{2\beta}u}{\partial x^{2\beta}}-\eta \frac{\partial^{2\beta}(u^2)}{\partial x^{2\beta}}-\zeta \frac{\partial^{4\beta}(u^2)}{\partial x^{4\beta}}-\mu \frac{\partial^{6\beta}(u^2)}{\partial x^{6\beta}}, ~t>0,~ \alpha, \beta \in (0,1],
\end{eqnarray}
where $\eta=4[\zeta-4 \mu], \zeta$ and $\mu$ are constants,
along with the initial conditions
\begin{eqnarray}\label{dbi}
u(x,0)=a+b \sin_{\beta}( x^{\beta})+c \cos_{\beta}( x^{\beta}), u_t(x,0)=0, ~~ a,b,c \in \mathbb{R}.
\end{eqnarray}
\end{ex}
Taking the Sumudu transform of both sides of Eq. \eqref{db}, we get
\begin{eqnarray}
S\Big[\frac{\partial^{2\alpha}u}{\partial t^{2\alpha}}\Big]= S\Big[\frac{\partial^{2\beta}u}{\partial x^{2\beta}}-\eta \frac{\partial^{2\beta}(u^2)}{\partial x^{2\beta}}-\zeta \frac{\partial^{4\beta}(u)^2}{\partial x^{4\beta}}-\mu \frac{\partial^{6\beta}(u)^2}{\partial x^{6\beta}}\Big].
\end{eqnarray} 

Using the property of Sumudu transform, we get
\begin{eqnarray}\label{dbs}
S[u(x,t)] &=& u(x,0)\nonumber \\ 
&+& \omega^{2\alpha}\Big(S\Big[\frac{\partial^{2\beta}u}{\partial x^{2\beta}}-\eta \frac{\partial^{2\beta}(u^2)}{\partial x^{2\beta}}-\zeta \frac{\partial^{4\beta}(u)^2}{\partial x^{4\beta}}-\mu \frac{\partial^{6\beta}(u)^2}{\partial x^{6\beta}}\Big] \Big).
\end{eqnarray}

Taking the inverse Sumudu transform of both sides of Eq. \eqref{dbs}
\begin{eqnarray}
u(x,t)&=& S^{-1}[u(x,0)]\nonumber \\ &+& S^{-1}\Big(\omega^{2\alpha}\Big(S\Big[\frac{\partial^{2\beta}u}{\partial x^{2\beta}}-\eta \frac{\partial^{2\beta}(u^2)}{\partial x^{2\beta}}-\zeta \frac{\partial^{4\beta}(u)^2}{\partial x^{4\beta}}-\mu \frac{\partial^{6\beta}(u)^2}{\partial x^{6\beta}}\Big]   \Big) \Big). \nonumber \\
\end{eqnarray}
Using the recurrence relation \eqref{nstim}, we get

\begin{eqnarray}
u_0= S^{-1}[u(x,0)]=a+b \sin_{\beta}( x^{\beta})+c \cos_{\beta}( x^{\beta}).
\end{eqnarray}
\begin{eqnarray}
u_1 &=& S^{-1}\Big(\omega^{2\alpha}\Big(S\Big[\frac{\partial^{2\beta}u_0}{\partial x^{2\beta}}-\eta \frac{\partial^{2\beta}(u_0^2)}{\partial x^{2\beta}}-\zeta \frac{\partial^{4\beta}(u_0)^2}{\partial x^{4\beta}}-\mu \frac{\partial^{6\beta}(u_0)^2}{\partial x^{6\beta}}\Big]   \Big) \Big)\nonumber \\
&=& \frac{t^{2 \alpha } (6 a (\zeta -5 \mu )-1) (b \cos_{\beta} ( x^{\beta})+c \sin_{\beta} ( x^{\beta}))}{\Gamma (2 \alpha +1)}.
\end{eqnarray}
Similarly,
\begin{eqnarray}
u_2 &=& \frac{t^{4 \alpha } (1-6 a (\zeta -5 \mu ))^2 (b \cos_{\beta} ( x^{\beta})+ c \sin_{\beta} ( x^{\beta})}{\Gamma (4 \alpha +1)}, \nonumber \\
u_3 &=& \frac{t^{6 \alpha } (6 a (\zeta -5 \mu )-1)^3 (b \cos_{\beta} ( x^{\beta})+c \sin_{\beta} ( x^{\beta}))}{\Gamma (6 \alpha +1)} \nonumber \\
\vdots
\end{eqnarray}
Hence, the series solution of Eqns (\ref{db}-\ref{dbi}) is given by
\begin{eqnarray}
u(x,t)  &=& u_0+u_1+u_2+..., \nonumber \\
&=& a+b \sin_{\beta}( x^{\beta})+c \cos_{\beta}(\lambda x^{\beta})\nonumber \\
&+& \frac{t^{2 \alpha } (6 a (\zeta -5 \mu )-1) (b \cos_{\beta} ( x^{\beta})+c \sin_{\beta} ( x^{\beta}))}{\Gamma (2 \alpha +1)} \nonumber \\
&+& \frac{t^{4 \alpha } (1-6 a (\zeta -5 \mu ))^2 (b \cos_{\beta} ( x^{\beta})+ c \sin_{\beta} ( x^{\beta})}{\Gamma (4 \alpha +1)} \nonumber \\
&+& \frac{t^{6 \alpha } (6 a (\zeta -5 \mu )-1)^3 (b \cos_{\beta} ( x^{\beta})+c \sin_{\beta} ( x^{\beta}))}{\Gamma (6 \alpha +1)} + ...,
\end{eqnarray}
which is equivalent to the following closed form solution:
\begin{eqnarray}
u(x,t) = a + b \sin_{\beta}( x^{\beta}) E_{2\alpha}(\delta t^{2\alpha})+c \cos_{\beta}( x^{\beta})E_{2\alpha}(\delta t^{2\alpha}), 
\end{eqnarray}
where $\delta= (6 a (\zeta -5 \mu )-1)$.

\begin{ex}
Consider the following general time space fractional diffusion-convection equation
\begin{eqnarray}\label{dce}
\frac{\partial^{\alpha}u}{\partial t^{\alpha}}=\Big(\frac{\partial ^{\beta}u}{\partial x^{\beta}} \Big)^2 \Big(\frac{\partial f(u)}{\partial u} \Big)+f(u) \frac{\partial^{2\beta}u}{\partial x^{2\beta}}-\frac{\partial^{\beta}u}{\partial x^{\beta}}\Big(\frac{\partial g(u)}{\partial u} \Big),~~ t>0, \alpha, \beta \in (0,1],
\end{eqnarray}
\end{ex}
where $f,g$ are the functions of $u$.
Here we consider some particular cases\\
Case 1: Let $f(u)=u, ~g(u)=k_1 = $ constant then Eq. \eqref{dce} reduces to 
\begin{eqnarray}\label{dce1}
\frac{\partial^{\alpha}u}{\partial t^{\alpha}}=\Big(\frac{\partial ^{\beta}u}{\partial x^{\beta}} \Big)^2 +u \frac{\partial^{2\beta}u}{\partial x^{2\beta}},
\end{eqnarray} 
along with the initial condition
\begin{eqnarray}\label{dce1i}
u(x,0)=a+bx^{\beta}.
\end{eqnarray}
Taking the Sumudu transform of both sides of Eq. \eqref{dce1}, we get
\begin{eqnarray}
S\Big[\frac{\partial^{\alpha}u}{\partial t^{\alpha}}\Big]=S\Big[\Big(\frac{\partial ^{\beta}u}{\partial x^{\beta}} \Big)^2 +u \frac{\partial^{2\beta}u}{\partial x^{2\beta}}\Big].
\end{eqnarray} 
Using the property of Sumudu transform
\begin{eqnarray}\label{dce1s}
S[u(x,t)]=u(x,0)+\omega^{\alpha}S\Big[\Big(\frac{\partial ^{\beta}u}{\partial x^{\beta}} \Big)^2 +u \frac{\partial^{2\beta}u}{\partial x^{2\beta}}\Big].
\end{eqnarray}
Taking inverse Sumudu tranform of both sides of Eq. \eqref{dce1s}

\begin{eqnarray}
u(x,t) = S^{-1}[u(x,0)]+S^{-1}\Big(\omega^{\alpha}S\Big[\Big(\frac{\partial ^{\beta}u}{\partial x^{\beta}} \Big)^2 +u \frac{\partial^{2\beta}u}{\partial x^{2\beta}}\Big]\Big).
\end{eqnarray}

Using the recurrence relation \eqref{nstim}, we get 
\begin{eqnarray}
u_0 &=& S^{-1}[u(x,0)]=a+bx^{\beta}, \nonumber \\
u_1 &=& S^{-1}\Big(\omega^{\alpha}S\Big[\Big(\frac{\partial ^{\beta}u_0}{\partial x^{\beta}} \Big)^2 +u_0 \frac{\partial^{2\beta}u_0}{\partial x^{2\beta}}\Big]\Big)=\frac{b^2 \Gamma (\beta +1)^2 t^{\alpha }}{\Gamma (\alpha +1)}, \nonumber \\
u_i &=& 0 ~~\forall~~ i \geq 2.
\end{eqnarray}
Hence, the exact solution of (\ref{dce1}-\ref{dce1i}) is given by 
\begin{eqnarray}
u(x,t) &=& a+\frac{b^2 \Gamma (\beta +1)^2 t^{\alpha }}{\Gamma (\alpha +1)}+b x^{\beta }.
\end{eqnarray} 

Case 2: Let $f(u)=\eta u$ and $g(u)=\frac{\zeta}{2}u^2$, where $\eta$ and $\zeta$ are constants and $\eta=\frac{\zeta}{2}$ then Eq. \eqref{dce} reduces to 
\begin{eqnarray}\label{dce2}
\frac{\partial^{\alpha}u}{\partial t^{\alpha}}=\eta \Big(\frac{\partial^{\beta}u}{\partial x^{\beta}}\Big)^2+ \eta u \frac{\partial^{2 \beta}u}{\partial x^{2 \beta}} -\zeta u \frac{\partial^{\beta}u}{\partial x^{\beta}},
\end{eqnarray} 
along with the initial condition
\begin{eqnarray}\label{dce2i}
u(x,0)=a+b E_{\beta}(x^{\beta}), ~ a,b \in \mathbb{R}
\end{eqnarray}

Taking the Sumudu transform of both sides of Eq.\eqref{dce2}
\begin{eqnarray}
S\Big[\frac{\partial^{\alpha}u}{\partial t^{\alpha}}\Big]=S\Big[\eta \Big(\frac{\partial^{\beta}u}{\partial x^{\beta}}\Big)^2+ \eta u \frac{\partial^{2 \beta}u}{\partial x^{2 \beta}} -\zeta u \frac{\partial^{\beta}u}{\partial x^{\beta}}\Big].
\end{eqnarray}
Using the property of Sumudu transform, we get
\begin{eqnarray}\label{dce2s}
S[u(x,t)]=u(x,0)+\omega^{\alpha}S\Big[\eta \Big(\frac{\partial^{\beta}u}{\partial x^{\beta}}\Big)^2+ \eta u \frac{\partial^{2 \beta}u}{\partial x^{2 \beta}} -\zeta u \frac{\partial^{\beta}u}{\partial x^{\beta}}\Big].
\end{eqnarray}
Taking inverse Sumudu transform of both sides of Eq.\eqref{dce2s}
\begin{eqnarray}
u(x,t)=S^{-1}[u(x,0)]+S^{-1}\Big(\omega^{\alpha}S\Big[\eta \Big(\frac{\partial^{\beta}u}{\partial x^{\beta}}\Big)^2+ \eta u \frac{\partial^{2 \beta}u}{\partial x^{2 \beta}} -\zeta u \frac{\partial^{\beta}u}{\partial x^{\beta}}\Big] \Big).
\end{eqnarray}
Using the recurrence relation \eqref{nstim}, we get
\begin{eqnarray}
u_0 &=& S^{-1}[u(x,0)] =a+b E_{\beta}(x^{\beta}), \nonumber \\
u_1 &=& S^{-1}\Big(\omega^{\alpha}S\Big[\eta \Big(\frac{\partial^{\beta}u_0}{\partial x^{\beta}}\Big)^2+ \eta u_0 \frac{\partial^{2 \beta}u}{\partial x^{2 \beta}} -\zeta u_0 \frac{\partial^{\beta}u}{\partial x^{\beta}}\Big] \Big), \nonumber \\
&=& -\frac{E_{\beta}(x^{\beta}) a b \zeta  t^{\alpha }}{2 \Gamma (\alpha +1)}, \nonumber \\
u_2 &=& \frac{E_{\beta}(x^{\beta}) a^2 b \zeta ^2 t^{2 \alpha }}{4 \Gamma (2 \alpha +1)}, \nonumber \\
u_3 &=& -\frac{E_{\beta}(x^{\beta}) a^3 b \zeta ^3 t^{3 \alpha }}{8 \Gamma (3 \alpha +1)},\nonumber \\
\vdots 
\end{eqnarray}
Hence, the series solution of Eqns. (\ref{dce2}-\ref{dce2i}) is given by
\begin{eqnarray}
u(x,t) &=& a+b E_{\beta}(x^{\beta}) -\frac{E_{\beta}(x^{\beta}) a b \zeta  t^{\alpha }}{2 \Gamma (\alpha +1)}+\frac{E_{\beta}(x^{\beta}) a^2 b \zeta ^2 t^{2 \alpha }}{4 \Gamma (2 \alpha +1)}\nonumber \\ &-&\frac{E_{\beta}(x^{\beta}) a^3 b \zeta ^3 t^{3 \alpha }}{8 \Gamma (3 \alpha +1)}+..., 
\end{eqnarray}
which is equivalent to the following closed from
\begin{eqnarray}
u(x,t)=a+b E_{\beta}(x^{\beta}) E_{\alpha}(-a \frac{\zeta}{2} t^{\alpha}).
\end{eqnarray}

\section{STIM for system of time and space FPDEs }\label{stimsys}
In this section we extend STIM to solve system of time and space fractional PDEs.

Consider the following system of time and space FPDEs:
\begin{equation}\label{ssteq}
\left.
\begin{aligned}
 \frac{\partial^{\gamma_i}u_i}{\partial t^{\gamma_i}} &= \mathscr{G}_i\Big(x, \bar{u},\frac{\partial^{\beta}\bar{u}}{\partial x^{\beta}},...,\frac{\partial^{l \beta}\bar{u}}{\partial x^{l \beta} } \Big),~  m_i-1<\gamma_i \leq m_i, \quad\\ 
     &~ ~i=1,2,...,q,~ n-1<\beta \leq n,~ m_i,l, n, q \in \mathbb{N},
\end{aligned}
\right\}
\end{equation}

along with the initial conditions
\begin{eqnarray}\label{ssini}
\frac{\partial^{j} u_{i}(x,0)}{\partial t^{j}}=g_{ij}(x), ~~j=0,1,2,...,m_i-1,
\end{eqnarray}
where  $\bar{u}=(u_1,u_2,...,u_q)$ and $\mathscr{G}_i\Big(x, \bar{u},\frac{\partial^{\beta}\bar{u}}{\partial x^{\beta}},...,\frac{\partial^{l \beta}\bar{u}}{\partial x^{l \beta} } \Big) $ is a linear/nonlinear operator.\\
After taking the Sumudu transform of both sides of Eq. \eqref{ssteq} and using Eq. \eqref{ssini}, we get
\begin{eqnarray}\label{ssumudu}
S[u_i(x,t)]=\sum_{j=0}^{m_i-1} \Big[\omega^{j}g_{ij}(x)\Big]+\omega^{\gamma_i}S\Big[\mathscr{G}_i\Big(x, \bar{u},\frac{\partial^{\beta}\bar{u}}{\partial x^{\beta}},...,\frac{\partial^{l \beta}\bar{u}}{\partial x^{l \beta} } \Big)\Big].
\end{eqnarray}
The inverse Sumudu transform of Eq. \eqref{ssumudu} yields the following system of equations
\begin{eqnarray}\label{sstn}
u_i(x,t) &=& S^{-1}\Big(\sum_{j=0}^{m_i-1} \Big[\omega^{j}g_{ij}(x)\Big]\Big) \nonumber \\
&+& S^{-1} \Big[\omega^{\gamma_i}S\Big(\mathscr{G}_i\Big(x, \bar{u},\frac{\partial^{\beta}\bar{u}}{\partial x^{\beta}},...,\frac{\partial^{l \beta}\bar{u}}{\partial x^{l \beta} } \Big)\Big)\Big], i=1,2,...,q.
\end{eqnarray}
Eq. \eqref{sstn} is of the following form
\begin{eqnarray}\label{sufn}
u_i(x,t) = f_i(x,t)+M_i\Big(x, \bar{u},\frac{\partial^{\beta}\bar{u}}{\partial x^{\beta}},...,\frac{\partial^{l \beta}\bar{u}}{\partial x^{l \beta}}\Big),
\end{eqnarray}
where 
\begin{equation}\label{asume}
\left.
\begin{aligned}
    f_i(x,t) &= S^{-1}\Big(\sum_{j=0}^{m_i-1}\Big[ \omega^{j}g_{ij}(x)\Big]\Big), \quad\\ 
   M_i\Big(x, \bar{u},\frac{\partial^{\beta}\bar{u}}{\partial x^{\beta}},...,\frac{\partial^{l \beta}\bar{u}}{\partial x^{l \beta}}\Big) &= S^{-1} \Big[\omega^{\gamma_i}S\Big(\mathscr{G}_i\Big(x, \bar{u},\frac{\partial^{\beta}\bar{u}}{\partial x^{\beta}},...,\frac{\partial^{l \beta}\bar{u}}{\partial x^{l \beta} } \Big)\Big)\Big].
\end{aligned}
\right\}
\end{equation}
Here $f_i$ is known function and $M_i$ is a linear/nonlinear operator. 
Functional equations of the form \eqref{sufn} can be solved by the DGJ decomposition method introduced by Daftardar-Gejji and Jafari \cite{daftardar2006iterative}.
DGJM represents the solution as an infinite series:
\begin{equation}\label{ssolution}
u_i=\sum_{j=0}^{\infty}u_i^{(j)}, ~~~ 1 \leq i \leq q,
\end{equation}
where the terms $u_i^{(j)}$ are calculated recursively. \\
Note that: Hence forward we use the following abbreviations:
\begin{eqnarray}
\bar{u}^{(j)}&=& (u_{1}^{(j)},u_{2}^{(j)},...,u_{q}^{(j)}), \nonumber \\
\sum_{j=0}^{r}\bar{u}^{(j)} &=& \Big(\sum_{j=0}^{r}u_1^{(j)},\sum_{j=0}^{r}u_2^{(j)},..., \sum_{j=0}^{r}u_q^{(j)} \Big),~ r \in \mathbb{N}\cup \{\infty\}, \nonumber \\
\frac{\partial^{k \beta}(\sum_{j=0}^{r}\bar{u}^{(j)})}{\partial x^{k \beta}}&=& \Big(\frac{\partial^{k\beta}(\sum_{j=0}^{r}u_1^{(j)})}{\partial x^{k\beta}}, \frac{\partial^{k\beta}(\sum_{j=0}^{r}u_2^{(j)})}{\partial x^{k\beta}},...,\frac{\partial^{k\beta}(\sum_{j=0}^{r}u_q^{(j)})}{\partial x^{k\beta}}\Big), k \in \mathbb{N}.\nonumber \\
\end{eqnarray}
The operator $M_i$ can be decomposed as:
\begin{eqnarray}
M_i\Big(x, \sum_{j=0}^{\infty}\bar{u}^{(j)},\frac{\partial^{\beta}(\sum_{j=0}^{\infty}\bar{u}^{(j)})}{\partial x^{\beta}},...,\frac{\partial^{l \beta}(\sum_{j=0}^{\infty}\bar{u}^{(j)})}{\partial x^{l \beta}}\Big)= \nonumber \\ M_i\Big(x,\bar{u}^{(0)},\frac{\partial^{\beta}\bar{u}^{(0)}}{\partial x^{\beta}},...,\frac{\partial^{l \beta}\bar{u}^{(0)}}{\partial x^{l \beta}}\Big) \nonumber \\ 
+\sum_{p=1}^\infty \Big(M_i\Big(x, \sum_{j=0}^{p}\bar{u}^{(j)},\frac{\partial^{\beta}(\sum_{j=0}^{p}\bar{u}^{(j)})}{\partial x^{\beta}},...,\frac{\partial^{l \beta}(\sum_{j=0}^{p}\bar{u}^{(j)})}{\partial x^{l \beta}}\Big)\Big)- \nonumber \\ \sum_{p=1}^\infty\Big(M_i\Big(x, \sum_{j=0}^{p-1}\bar{u}^{(j)},\frac{\partial^{\beta}(\sum_{j=0}^{p-1}\bar{u}^{(j)})}{\partial x^{\beta}},...,\frac{\partial^{l \beta}(\sum_{j=0}^{p-1}\bar{u}^{(j)})}{\partial x^{l \beta}}\Big)\Big).
\end{eqnarray}
Therefore,
\begin{equation}\label{sdecom}
\left.
\begin{aligned}
   S^{-1} \Big[\omega^{\gamma_i}S\Big(\mathscr{G}_i\Big(x,\sum_{j=0}^{\infty}\bar{u}^{(j)} ,\frac{\partial^{\beta}(\sum_{j=0}^{\infty} \bar{u}^{(j)})}{\partial x^{\beta}}, ...,\frac{\partial^{l \beta}(\sum_{j=0}^{\infty} \bar{u}^{(j)})}{\partial x^{l \beta} } \Big)\Big)\Big] & \quad\\ 
     = S^{-1} \Big[\omega^{\gamma_i}S\Big(\mathscr{G}_i\Big(x, \bar{u}^{(0)},\frac{\partial^{\beta}\bar{u}^{(0)}}{\partial x^{\beta}}, ...,\frac{\partial^{l \beta}\bar{u}^{(0)}}{\partial x^{l \beta} }\Big)\Big)\Big]& \quad \\
    +\sum_{p=1}^{\infty} S^{-1} \Big[\omega^{\gamma_i}S\Big(\mathscr{G}_i\Big(x,\sum_{j=0}^{p}\bar{u}^{(j)} ,\frac{\partial^{\beta}(\sum_{j=0}^{p} \bar{u}^{(j)})}{\partial x^{\beta}}, ...,\frac{\partial^{l \beta}(\sum_{j=0}^{p} \bar{u}^{(j)})}{\partial x^{l \beta} }\Big)\Big)\Big]&\quad \\
     -\sum_{p=1}^{\infty} S^{-1} \Big[\omega^{\gamma_i}S\Big(\mathscr{G}_i\Big(x, \sum_{j=0}^{p-1}\bar{u}^{(j)} ,\frac{\partial^{\beta}(\sum_{j=0}^{p-1} \bar{u}^{(j)})}{\partial x^{\beta}},...,\frac{\partial^{l \beta}(\sum_{j=0}^{p-1}\bar{u}^{(j)})}{\partial x^{l \beta} }\Big)\Big)\Big].
\end{aligned}
\right\}
\end{equation}

Using Eqns. (\ref{ssolution},\ref{sdecom}) in Eq. \eqref{sufn}, we get
\begin{eqnarray}
\sum_{j=0}^{\infty}u_i^{(j)} = S^{-1}\Big(\sum_{j=0}^{m_i-1} \Big[\omega^{j}g_{ij}(x)\Big]\Big)+S^{-1} \Big[\omega^{\gamma_i}S\Big(\mathscr{G}_i\Big(x, \bar{u}^{(0)},\frac{\partial^{\beta}\bar{u}^{(0)}}{\partial x^{\beta}}, ...,\frac{\partial^{l \beta}\bar{u}^{(0)}}{\partial x^{l \beta} }\Big)\Big)\Big] \nonumber \\
+\sum_{p=1}^{\infty} \Big(S^{-1} \Big[\omega^{\gamma_i}S\Big(\mathscr{G}_i\Big(x, \sum_{j=0}^{p}\bar{u}^{(j)} ,\frac{\partial^{\beta}(\sum_{j=0}^{p} \bar{u}^{(j)})}{\partial x^{\beta}},...,\frac{\partial^{l \beta}(\sum_{j=0}^{p}\bar{u}^{(j)})}{\partial x^{l \beta} }\Big)\Big)\Big] \nonumber \\
- S^{-1} \Big[\omega^{\gamma_i}S\Big(\mathscr{G}_i\Big(x, \sum_{j=0}^{p-1}\bar{u}^{(j)} ,\frac{\partial^{\beta}(\sum_{j=0}^{p-1} \bar{u}^{(j)})}{\partial x^{\beta}}, ...,\frac{\partial^{l \beta}(\sum_{j=0}^{p-1} \bar{u}^{(j)})}{\partial x^{l \beta} }\Big)\Big)\Big]\Big).
\end{eqnarray}
We define the recurrence relation as follows:

\begin{equation}\label{nsstim}
\left.
\begin{aligned}
u_i^{(0)} &= S^{-1}\Big(\sum_{j=0}^{m_i-1} \Big[\omega^{j}g_{ij}(x)\Big]\Big),\\ 
u_i^{(1)} &= S^{-1} \Big[\omega^{\gamma_i}S\Big(\mathscr{G}_i\Big(x, \bar{u}^{(0)},\frac{\partial^{\beta}\bar{u}^{(0)}}{\partial x^{\beta}}, ...,\frac{\partial^{l \beta}\bar{u}^{(0)}}{\partial x^{l \beta} }\Big)\Big)\Big],\\
\hspace*{-1.5cm}u_i^{(m+1)}&= 
S^{-1} \Big[\omega^{\gamma_i}S\Big(\mathscr{G}_i\Big(x, \sum_{j=0}^{m}\bar{u}^{(j)} ,\frac{\partial^{\beta}(\sum_{j=0}^{m} \bar{u}^{(j)})}{\partial x^{\beta}}, ...,\frac{\partial^{l \beta}(\sum_{j=0}^{m} \bar{u}^{(j)})}{\partial x^{l \beta} }\Big)\Big)\Big]\\
&- S^{-1} \Big[\omega^{\gamma_i}S\Big(\mathscr{G}_i\Big(x, \sum_{j=0}^{m-1}\bar{u}^{(j)} ,\frac{\partial^{\beta}(\sum_{j=0}^{m-1} \bar{u}^{(j)})}{\partial x^{\beta}}, ...,\frac{\partial^{l \beta}(\sum_{j=0}^{m-1} \bar{u}^{(j)})}{\partial x^{l \beta} }\Big)\Big)\Big],\\&~~~~~~~~~~for ~~~~ m\geq 1.
\end{aligned}
\right\}
\end{equation}

The m-term approximate solution of Eqns. (\ref{ssteq}-\ref{ssini}) is given by $u_i \approx u_i^{(0)}+u_i^{(1)}+ \cdots +u_i^{(m-1)}$ or $u_i \approx u_{i0}+u_{i1}+ \cdots +u_{i(m-1)}$.

\section{Illustrative Examples}\label{stimsysexample}
In this section we solve system of time and space FPDEs using STIM  derived in \ref{stimsys}.
\begin{ex}
Consider the following system of time and space fractional Boussinesq PDEs ($t>0, 0<\alpha_1, \alpha_2, \beta \leq 1$):  
\begin{eqnarray}\label{bous}
\frac{\partial^{\alpha_1} u_1}{\partial t^{\alpha_1}}&=& -\frac{\partial^\beta u_2}{\partial x^{\beta}}, \nonumber \\
\frac{\partial^{\alpha_2}u_2}{\partial t^{\alpha_2}}&=& -m_1\frac{\partial^{\beta}u_1}{\partial x^{\beta}}+3u_1\Big(\frac{\partial^{\beta}u_1}{\partial x^{\beta}}\Big)+m_2\frac{\partial^{3\beta}u_1}{\partial x^{3\beta}},
\end{eqnarray}
along with the following initial conditions
\begin{eqnarray}\label{bousi}
u_1(x,0)=a+ b x^{\beta}, u_2(x,0)=c,~~ a,b,c \in \mathbb{R}.
\end{eqnarray}
\end{ex}
Taking the Sumudu transform of both sides of Eqns.\eqref{bous}
\begin{eqnarray}
S \Big[\frac{\partial^{\alpha_1} u_1}{\partial t^{\alpha_1}}\Big]&=& S\Big[-\frac{\partial^\beta u_2}{\partial x^{\beta}}\Big],\nonumber \\
S\Big[\frac{\partial^{\alpha_2}u_2}{\partial t^{\alpha_2}}\Big]&=& S\Big[-m_1\frac{\partial^{\beta}u_1}{\partial x^{\beta}}+3u_1\Big(\frac{\partial^{\beta}u_1}{\partial x^{\beta}}\Big)+m_2\frac{\partial^{3\beta}u_1}{\partial x^{3\beta}}\Big].
\end{eqnarray}
In view of \eqref{sproperty}, we get

\begin{eqnarray}\label{bouss}
S[u_1(x,t)]&=& u_1(x,0)+ \omega^{\alpha_1} S\Big[-\frac{\partial^\beta u_2}{\partial x^{\beta}}\Big],\nonumber \\
S[u_2(x,t)]&=& u_2(x,0)+ \omega^{\alpha_2}S\Big[-m_1\frac{\partial^{\beta}u_1}{\partial x^{\beta}}+3u_1\Big(\frac{\partial^{\beta}u_1}{\partial x^{\beta}}\Big)+m_2\frac{\partial^{3\beta}u_1}{\partial x^{3\beta}}\Big].\nonumber \\
\end{eqnarray}

Taking the inverse Sumudu transform of both sides of Eqns. \eqref{bouss}
 \begin{eqnarray}
u_1(x,t)&=& S^{-1}[u_1(x,0)]+ S^{-1}\Big(\omega^{\alpha_1} S\Big[-\frac{\partial^\beta u_2}{\partial x^{\beta}}\Big]\Big),\nonumber \\
u_2(x,t)&=& S^{-1}[u_2(x,0)]\nonumber \\ &+& S^{-1}\Big( \omega^{\alpha_2}S\Big[-m_1\frac{\partial^{\beta}u_1}{\partial x^{\beta}}+3u_1\Big(\frac{\partial^{\beta}u_1}{\partial x^{\beta}}\Big)+m_2\frac{\partial^{3\beta}u_1}{\partial x^{3\beta}}\Big]\Big).
 \end{eqnarray}
The recurrence relation \eqref{nsstim} yields
\begin{eqnarray}
u_{10} &=& S^{-1}[u_1(x,0)]=a+ b x^{\beta},\nonumber\\
u_{20} &=& S^{-1}[u_2(x,0)]=c, \nonumber \\
u_{11} &=& S^{-1}\Big(\omega^{\alpha_1} S\Big[-\frac{\partial^\beta u_{20}}{\partial x^{\beta}}\Big]\Big)=0,\nonumber \\
u_{21} &=& S^{-1}\Big( \omega^{\alpha_2}S\Big[-m_1\frac{\partial^{\beta}u_{10}}{\partial x^{\beta}}+3u_{10}\Big(\frac{\partial^{\beta}u_{10}}{\partial x^{\beta}}\Big)+m_2\frac{\partial^{3\beta}u_{10}}{\partial x^{3\beta}}\Big]\Big), \nonumber \\
&=& \frac{3 a b \Gamma (\beta +1) t^{\alpha _2}}{\Gamma \left(\alpha _2+1\right)}+\frac{3 b^2 \Gamma (\beta +1) t^{\alpha _2} x^{\beta }}{\Gamma \left(\alpha _2+1\right)}-\frac{b m_1 \Gamma (\beta +1) t^{\alpha _2}}{\Gamma \left(\alpha _2+1\right)}, \nonumber \\
u_{12} &=& -\frac{3 b^2 \Gamma (\beta +1)^2 t^{\alpha _1+\alpha _2}}{\Gamma \left(\alpha _1+\alpha _2+1\right)}, \nonumber \\
u_{22} &=& 0, \nonumber \\
u_{13} &=& 0, \nonumber \\
u_{23} &=& -\frac{9 b^3 \Gamma (\beta +1)^3 t^{\alpha _1+2 \alpha _2}}{\Gamma \left(\alpha _1+2 \alpha _2+1\right)}, \nonumber \\
u_{1n} &=& 0, n \geq 4, \nonumber \\
u_{2n} &=& 0, n \geq 4.
\end{eqnarray}
Hence, the exact solution of the system(\ref{bous}) along with the initial conditions \eqref{bousi} is given by 
\begin{eqnarray}\label{boussol}
u_1(x,t)&=& u_{10}+u_{11}+u_{12}+u_{13}, \nonumber \\
&=& a-\frac{3 b^2 \Gamma (\beta +1)^2 t^{\alpha _1+\alpha _2}}{\Gamma \left(\alpha _1+\alpha _2+1\right)}+b x^{\beta }, \nonumber \\
u_2(x,t)&=& u_{20}+u_{21}+u_{22}+u_{23}, \nonumber \\
&=& c+ \frac{3 a b \Gamma (\beta +1) t^{\alpha _2}}{\Gamma \left(\alpha _2+1\right)}-\frac{9 b^3 \Gamma (\beta +1)^3 t^{\alpha _1+2 \alpha _2}}{\Gamma \left(\alpha _1+2 \alpha _2+1\right)}\nonumber \\ 
&+&\frac{3 b^2 \Gamma (\beta +1) t^{\alpha _2} x^{\beta }}{\Gamma \left(\alpha _2+1\right)}-\frac{b m_1 \Gamma (\beta +1) t^{\alpha _2}}{\Gamma \left(\alpha _2+1\right)}.
\end{eqnarray}
In case of when $a=e, b=2$ and $c=3/2$, this solution is same as obtained using invariant subspace method in \cite{choudhary2018solving}. 

\begin{ex}
Consider the following two-coupled time and space fractional diffusion system
\begin{eqnarray}\label{diffusion}
\frac{\partial^{\alpha_1}u_1}{\partial t^{\alpha_1}}&=&\frac{\partial^{2 \beta}u_1}{\partial x^{2\beta}}+\mu \frac{\partial^{\beta}}{\partial x^{\beta}}\Big(u_2 \frac{\partial^{\beta}u_2}{\partial x^{\beta}}\Big)+\xi u_2^2, \nonumber \\
\frac{\partial^{\alpha_2}u_2}{\partial t^{\alpha_2}}&=&\frac{\partial^{2\beta}u_2}{\partial x^{2\beta}}+\eta \frac{\partial^{2 \beta}u_1}{\partial x^{2\beta}}+\zeta u_1 +\delta u_2, t>0, 0<\alpha_1,\alpha_2,\beta \leq 1,
\end{eqnarray} 
where $\mu, \xi, \eta, \zeta, \delta$ are arbitrary constants, $\mu$ and $\xi$ are not simultaneously zero, we consider $\xi=-2 \mu \lambda^2, \zeta =\eta \kappa^2$, along with the initial conditions
\begin{eqnarray}\label{diffusioni}
\hspace{-1cm}u_1(x,0)=a \cos_{\beta} (\kappa  x^{\beta})+b \sin_{\beta} (\kappa  x^{\beta}), u_2(x,0)=c E_{\beta}(-\lambda x^{\beta}),~ a,b,c,\lambda,\kappa \in \mathbb{R}.
\end{eqnarray} 
\end{ex}
Taking the Sumudu transform on both sides of Eqns.\eqref{diffusion}
\begin{eqnarray}
S\Big[\frac{\partial^{\alpha_1}u_1}{\partial t^{\alpha_1}}\Big]&=& S\Big[\frac{\partial^{2 \beta}u_1}{\partial x^{2\beta}}+\mu \frac{\partial^{\beta}}{\partial x^{\beta}}\Big(u_2 \frac{\partial^{\beta}u_2}{\partial x^{\beta}}\Big)+\xi u_2^2\Big], \nonumber \\
S\Big[\frac{\partial^{\alpha_2}u_2}{\partial t^{\alpha_2}}\Big]&=& S\Big[\frac{\partial^{2\beta}u_2}{\partial x^{2\beta}}+\eta \frac{\partial^{2 \beta}u_1}{\partial x^{2\beta}}+\zeta u_1 +\delta u_2\Big].
\end{eqnarray}
After using the property of Sumudu transform \eqref{sproperty} we get,
\begin{eqnarray}
S[u_1(x,t)]&=& u_1(x,0)+\omega^{\alpha_1} S\Big[\frac{\partial^{2 \beta}u_1}{\partial x^{2\beta}}+\mu \frac{\partial^{\beta}}{\partial x^{\beta}}\Big(u_2 \frac{\partial^{\beta}u_2}{\partial x^{\beta}}\Big)+\xi u_2^2\Big], \nonumber \\
S[u_2(x,t)]&=& u_2(x,0)+\omega^{\alpha_2} S\Big[\frac{\partial^{2\beta}u_2}{\partial x^{2\beta}}+\eta \frac{\partial^{2 \beta}u_1}{\partial x^{2\beta}}+\zeta u_1 +\delta u_2\Big].
\end{eqnarray}
Taking the inverse Sumudu transform 
\begin{eqnarray}
u_1(x,t)&=& S^{-1}[u_1(x,0)]+S^{-1}\Big(\omega^{\alpha_1} S\Big[\frac{\partial^{2 \beta}u_1}{\partial x^{2\beta}}+\mu \frac{\partial^{\beta}}{\partial x^{\beta}}\Big(u_2 \frac{\partial^{\beta}u_2}{\partial x^{\beta}}\Big)+\xi u_2^2\Big]\Big), \nonumber \\
u_2(x,t)&=& S^{-1}[u_2(x,0)]+S^{-1}\Big(\omega^{\alpha_2} S\Big[\frac{\partial^{2\beta}u_2}{\partial x^{2\beta}}+\eta \frac{\partial^{2 \beta}u_1}{\partial x^{2\beta}}+\zeta u_1 +\delta u_2\Big]\Big). \nonumber \\
\end{eqnarray}
Using the recurrence relation \eqref{nsstim}, we get
\begin{eqnarray}
u_{10}&=& S^{-1}[u_1(x,0)]=a \cos_{\beta} (\kappa  x^{\beta})+b \sin_{\beta} (\kappa  x^{\beta}), \nonumber \\
u_{20}&=& S^{-1}[u_2(x,0)]=c E_{\beta}(-\lambda x^{\beta}), \nonumber\\
u_{11}&=& S^{-1}\Big(\omega^{\alpha_1} S\Big[\frac{\partial^{2 \beta}u_{10}}{\partial x^{2\beta}}+\mu \frac{\partial^{\beta}}{\partial x^{\beta}}\Big(u_{20} \frac{\partial^{\beta}u_{20}}{\partial x^{\beta}}\Big)+\xi u_{20}^2\Big]\Big),\nonumber \\
&=& -\frac{a \kappa ^2 t^{\alpha _1} \cos_{\beta} (\kappa  x^{\beta})}{\Gamma \left(\alpha _1+1\right)}-\frac{b \kappa ^2 t^{\alpha _1} \sin_{\beta} (\kappa  x^{\beta})}{\Gamma \left(\alpha _1+1\right)},\nonumber \\
u_{21} &=& S^{-1}\Big(\omega^{\alpha_2} S\Big[\frac{\partial^{2\beta}u_{20}}{\partial x^{2\beta}}+\eta \frac{\partial^{2 \beta}u_{10}}{\partial x^{2\beta}}+\zeta u_{10} +\delta u_{20}\Big]\Big), \nonumber \\
&=& \frac{c \delta  t^{\alpha _2}E_{\beta}(-\lambda x^{\beta})}{\Gamma \left(\alpha _2+1\right)}+\frac{c \lambda ^2 t^{\alpha _2} E_{\beta}(-\lambda x^{\beta})}{\Gamma \left(\alpha _2+1\right)}, \nonumber \\
u_{12}&=& \frac{a \kappa ^4 t^{2 \alpha _1} \cos_{\beta} (\kappa  x^{\beta})}{\Gamma \left(2 \alpha _1+1\right)}+\frac{b \kappa ^4 t^{2 \alpha _1} \sin_{\beta} (\kappa  x^{\beta})}{\Gamma \left(2 \alpha _1+1\right)}, \nonumber \\
u_{22}&=& \frac{c \delta ^2 t^{2 \alpha _2} E_{\beta}(-\lambda x^{\beta})}{\Gamma \left(2 \alpha _2+1\right)}+\frac{2 c \delta  \lambda ^2 t^{2 \alpha _2} E_{\beta}(-\lambda x^{\beta})}{\Gamma \left(2 \alpha _2+1\right)}+\frac{c \lambda ^4 t^{2 \alpha _2} E_{\beta}(-\lambda x^{\beta})}{\Gamma \left(2 \alpha _2+1\right)}, \nonumber \\
u_{13}&=& -\frac{a \kappa ^6 t^{3 \alpha _1} \cos_{\beta} (\kappa  x^{\beta})}{\Gamma \left(3 \alpha _1+1\right)}-\frac{b \kappa ^6 t^{3 \alpha _1} \sin_{\beta} (\kappa  x^{\beta})}{\Gamma \left(3 \alpha _1+1\right)} ,\nonumber \\
\end{eqnarray}
\begin{eqnarray}
u_{23}&=& \frac{c \delta ^3 t^{3 \alpha _2} E_{\beta}(-\lambda x^{\beta})}{\Gamma \left(3 \alpha _2+1\right)}+\frac{3 c \delta ^2 \lambda ^2 t^{3 \alpha _2} E_{\beta}(-\lambda x^{\beta})}{\Gamma \left(3 \alpha _2+1\right)}+\frac{3 c \delta  \lambda ^4 t^{3 \alpha _2} E_{\beta}(-\lambda x^{\beta})}{\Gamma \left(3 \alpha _2+1\right)}\nonumber \\
&+&\frac{c \lambda ^6 t^{3 \alpha _2} E_{\beta}(-\lambda x^{\beta})}{\Gamma \left(3 \alpha _2+1\right)}, \nonumber \\
\vdots
\end{eqnarray}
Hence, the series solution of two-coupled time and space fractional diffusion system (\ref{diffusion}-\ref{diffusioni}) is given by
\begin{eqnarray}
u_1(x,t) &=& a \cos_{\beta} (\kappa  x^{\beta})+b \sin_{\beta} (\kappa  x^{\beta}) -\frac{a \kappa ^2 t^{\alpha _1} \cos_{\beta} (\kappa  x^{\beta})}{\Gamma \left(\alpha _1+1\right)}-\frac{b \kappa ^2 t^{\alpha _1} \sin_{\beta} (\kappa  x^{\beta})}{\Gamma \left(\alpha _1+1\right)} \nonumber \\ &+& \frac{a \kappa ^4 t^{2 \alpha _1} \cos_{\beta} (\kappa  x^{\beta})}{\Gamma \left(2 \alpha _1+1\right)}+\frac{b \kappa ^4 t^{2 \alpha _1} \sin_{\beta} (\kappa  x^{\beta})}{\Gamma \left(2 \alpha _1+1\right)}\nonumber \\&-&\frac{a \kappa ^6 t^{3 \alpha _1} \cos_{\beta} (\kappa  x^{\beta})}{\Gamma \left(3 \alpha _1+1\right)}-\frac{b \kappa ^6 t^{3 \alpha _1} \sin_{\beta} (\kappa  x^{\beta})}{\Gamma \left(3 \alpha _1+1\right)}+..., \\
u_2(x,t)&=& c E_{\beta}(-\lambda x^{\beta})+\frac{c \delta  t^{\alpha _2}E_{\beta}(-\lambda x^{\beta})}{\Gamma \left(\alpha _2+1\right)}+\frac{c \lambda ^2 t^{\alpha _2} E_{\beta}(-\lambda x^{\beta})}{\Gamma \left(\alpha _2+1\right)} \nonumber \\ 
&+&\frac{c \delta ^2 t^{2 \alpha _2} E_{\beta}(-\lambda x^{\beta})}{\Gamma \left(2 \alpha _2+1\right)}+\frac{2 c \delta  \lambda ^2 t^{2 \alpha _2} E_{\beta}(-\lambda x^{\beta})}{\Gamma \left(2 \alpha _2+1\right)}+\frac{c \lambda ^4 t^{2 \alpha _2} E_{\beta}(-\lambda x^{\beta})}{\Gamma \left(2 \alpha _2+1\right)}\nonumber \\ 
&+&\frac{c \delta ^3 t^{3 \alpha _2} E_{\beta}(-\lambda x^{\beta})}{\Gamma \left(3 \alpha _2+1\right)}+\frac{3 c \delta ^2 \lambda ^2 t^{3 \alpha _2} E_{\beta}(-\lambda x^{\beta})}{\Gamma \left(3 \alpha _2+1\right)}+\frac{3 c \delta  \lambda ^4 t^{3 \alpha _2} E_{\beta}(-\lambda x^{\beta})}{\Gamma \left(3 \alpha _2+1\right)}\nonumber \\
&+&\frac{c \lambda ^6 t^{3 \alpha _2} E_{\beta}(-\lambda x^{\beta})}{\Gamma \left(3 \alpha _2+1\right)}+\cdots.
\end{eqnarray}
These series solutions converge to:
\begin{eqnarray}
u_1(x,t)&=& [a \cos_{\beta} (\kappa  x^{\beta})+b \sin_{\beta} (\kappa  x^{\beta})] E_{\alpha_1}(-\kappa^2 t^{\alpha_1} ), \nonumber \\
u_2(x,t)&=& c E_{\alpha_2}[(\delta+\lambda^2)t^{\alpha_2}] E_{\beta}(-\lambda x^{\beta}).  
\end{eqnarray} 

\begin{ex}
Consider the following two-coupled time and space fractional PDE:
\begin{eqnarray}\label{pde1}
\frac{\partial^{\alpha}u_1}{\partial t^{\alpha}}&=& \frac{\partial^{\beta}}{\partial x^{\beta}}\Big(\frac{\partial^{4\beta}u_1}{\partial x^{4\beta}}+\eta u_2 \frac{\partial^{\beta}u_2}{\partial x^{\beta}} \Big)+\zeta u_2^2, \nonumber \\
\frac{\partial^{\alpha}u_2}{\partial t^{\alpha}}&=& \frac{\partial^{4\beta}u_1}{\partial x^{4\beta}}+\delta u_1+\tau u_2,~~t>0,0<\alpha,\beta \leq 1,
\end{eqnarray}
here $\eta, \zeta, \delta,\tau$ all are arbitrary constants, $\eta$ and  $\zeta$ are not simultaneously zero (taking $\zeta=-2 \eta$), along with the initial conditions
\begin{eqnarray}\label{pde1i}
u_1(x,0)=b E_{\beta}(-x^{\beta}),~ u_2(x,0)= d E_{\beta}(-x^{\beta}),~~ b,d \in \mathbb{R}.
\end{eqnarray}
\end{ex}
Taking the Sumudu transform of both sides of Eqns. \eqref{pde1}
\begin{eqnarray}
S\Big[\frac{\partial^{\alpha}u_1}{\partial t^{\alpha}}\Big]&=& S\Big[\frac{\partial^{\beta}}{\partial x^{\beta}}\Big(\frac{\partial^{4\beta}u_1}{\partial x^{4\beta}}+\eta u_2 \frac{\partial^{\beta}u_2}{\partial x^{\beta}} \Big)+\zeta u_2^2\Big], \nonumber \\
S\Big[\frac{\partial^{\alpha}u_2}{\partial t^{\alpha}}\Big]&=& S\Big[\frac{\partial^{4\beta}u_1}{\partial x^{4\beta}}+\delta u_1+\tau u_2\Big].
\end{eqnarray}

After using the property of Sumudu transform\eqref{sproperty}, we get
\begin{eqnarray}\label{pde1s}
S[u_1(x,t)]&=& u_1(x,0)+\omega^{\alpha}S\Big[\frac{\partial^{\beta}}{\partial x^{\beta}}\Big(\frac{\partial^{4\beta}u_1}{\partial x^{4\beta}}+\eta u_2 \frac{\partial^{\beta}u_2}{\partial x^{\beta}} \Big)+\zeta u_2^2\Big], \nonumber \\
S[u_2(x,t)]&=& u_2(x,0)+ \omega^{\alpha}S\Big[\frac{\partial^{4\beta}u_1}{\partial x^{4\beta}}+\delta u_1+\tau u_2\Big].
\end{eqnarray}
Taking the inverse Sumudu transform of both sides of Eqns. \eqref{pde1s}
\begin{eqnarray}
u_1(x,t)&=& S^{-1}[u_1(x,0)]+S^{-1}\Big(\omega^{\alpha}S\Big[\frac{\partial^{\beta}}{\partial x^{\beta}}\Big(\frac{\partial^{4\beta}u_1}{\partial x^{4\beta}}+\eta u_2 \frac{\partial^{\beta}u_2}{\partial x^{\beta}} \Big)+\zeta u_2^2\Big]\Big), \nonumber \\
u_2(x,t)&=& S^{-1}[u_2(x,0)]+ S^{-1}\Big(\omega^{\alpha}S\Big[\frac{\partial^{4\beta}u_1}{\partial x^{4\beta}}+\delta u_1+\tau u_2\Big]\Big).
\end{eqnarray}

In view of the recurrence relation \eqref{nsstim}

\begin{eqnarray}
u_{10}&=& S^{-1}[u_1(x,0)]=E_{\beta}(-x^{\beta}), \nonumber \\
u_{20}&=& S^{-1}[u_2(x,0)]= E_{\beta}(-x^{\beta}), \nonumber \\
u_{11}&=& S^{-1}\Big(\omega^{\alpha}S\Big[\frac{\partial^{\beta}}{\partial x^{\beta}}\Big(\frac{\partial^{4\beta}u_{10}}{\partial x^{4\beta}}+\eta u_{20} \frac{\partial^{\beta}u_{20}}{\partial x^{\beta}} \Big)+\zeta u_{20}^2\Big]\Big), \nonumber \\
&=&  -\frac{b E_{\beta}(-x^{\beta}) t^{\alpha }}{\Gamma (\alpha +1)}, \nonumber \\
u_{21}&=& S^{-1}\Big(\omega^{\alpha}S\Big[\frac{\partial^{4\beta}u_{10}}{\partial x^{4\beta}}+\delta u_{10}+\tau u_{20}\Big]\Big), \nonumber \\
&=& \frac{b \delta E_{\beta}(-x^{\beta}) t^{\alpha }}{\Gamma (\alpha +1)}+\frac{b E_{\beta}(-x^{\beta}) t^{\alpha }}{\Gamma (\alpha +1)}+\frac{d \tau E_{\beta}(-x^{\beta}) t^{\alpha }}{\Gamma (\alpha +1)}, \nonumber \\
u_{12}&=& \frac{b E_{\beta}(-x^{\beta}) t^{2 \alpha }}{\Gamma (2 \alpha +1)}, \nonumber \\
u_{22} &=& \frac{b \delta \tau E_{\beta}(-x^{\beta}) t^{2 \alpha }}{\Gamma (2 \alpha +1)}-\frac{b \delta E_{\beta}(-x^{\beta}) t^{2 \alpha }}{\Gamma (2 \alpha +1)}+\frac{b \tau E_{\beta}(-x^{\beta}) t^{2 \alpha }}{\Gamma (2 \alpha +1)}\nonumber \\ &-& \frac{b E_{\beta}(-x^{\beta}) t^{2 \alpha }}{\Gamma (2 \alpha +1)}+\frac{d \tau^2 E_{\beta}(-x^{\beta}) t^{2 \alpha }}{\Gamma (2 \alpha +1)}, \nonumber \\
\end{eqnarray}
\begin{eqnarray}
u_{13}&=&-\frac{b E_{\beta}(-x^{\beta}) t^{3 \alpha }}{\Gamma (3 \alpha +1)},
\nonumber \\
u_{23}&=& \frac{b \delta \tau^2 E_{\beta}(-x^{\beta}) t^{3 \alpha }}{\Gamma (3 \alpha +1)}-\frac{b \delta \tau E_{\beta}(-x^{\beta}) t^{3 \alpha }}{\Gamma (3 \alpha +1)}+\frac{b \delta E_{\beta}(-x^{\beta}) t^{3 \alpha }}{\Gamma (3 \alpha +1)}+\frac{b \tau^2 E_{\beta}(-x^{\beta}) t^{3 \alpha }}{\Gamma (3 \alpha +1)}\nonumber \\ &-& \frac{b \tau E_{\beta}(-x^{\beta}) t^{3 \alpha }}{\Gamma (3 \alpha +1)}+\frac{b E_{\beta}(-x^{\beta}) t^{3 \alpha }}{\Gamma (3 \alpha +1)}+\frac{d \tau^3 E_{\beta}(-x^{\beta}) t^{3 \alpha }}{\Gamma (3 \alpha +1)},\nonumber \\
u_{14}&=&\frac{b E_{\beta}(-x^{\beta}) t^{4 \alpha }}{\Gamma (4 \alpha +1)},\nonumber \\
u_{24}&=& \frac{b \delta \tau^3 E_{\beta}(-x^{\beta}) t^{4 \alpha }}{\Gamma (4 \alpha +1)}-\frac{b \delta \tau^2 E_{\beta}(-x^{\beta}) t^{4 \alpha }}{\Gamma (4 \alpha +1)}\nonumber +\frac{b \delta \tau E_{\beta}(-x^{\beta}) t^{4 \alpha }}{\Gamma (4 \alpha +1)}\nonumber \\ &-&\frac{b \delta E_{\beta}(-x^{\beta}) t^{4 \alpha }}{\Gamma (4 \alpha +1)}+ \frac{b \tau^3 E_{\beta}(-x^{\beta}) t^{4 \alpha }}{\Gamma (4 \alpha +1)}-\frac{b \tau^2 E_{\beta}(-x^{\beta}) t^{4 \alpha }}{\Gamma (4 \alpha +1)}+\frac{b \tau E_{\beta}(-x^{\beta}) t^{4 \alpha }}{\Gamma (4 \alpha +1)}\nonumber \\ &-&\frac{b E_{\beta}(-x^{\beta}) t^{4 \alpha }}{\Gamma (4 \alpha +1)}+\frac{d \tau^4 E_{\beta}(-x^{\beta}) t^{4 \alpha }}{\Gamma (4 \alpha +1)},\nonumber \\
\vdots
\end{eqnarray}
Hence, the series solution of time and space fractional PDE \eqref{pde1} along with the initial conditions \eqref{pde1i} is given by
\begin{eqnarray}
u_1(x,t)&=& E_{\beta}(-x^{\beta})-\frac{b E_{\beta}(-x^{\beta}) t^{\alpha }}{\Gamma (\alpha +1)}+\frac{b E_{\beta}(-x^{\beta}) t^{2 \alpha }}{\Gamma (2 \alpha +1)}-\frac{b E_{\beta}(-x^{\beta}) t^{3 \alpha }}{\Gamma (3 \alpha +1)}\nonumber \\ &+& \frac{b E_{\beta}(-x^{\beta}) t^{4 \alpha }}{\Gamma (4 \alpha +1)}-\cdots. \nonumber \\
u_2(x,t)&=&  E_{\beta}(-x^{\beta})+\frac{b \delta E_{\beta}(-x^{\beta}) t^{\alpha }}{\Gamma (\alpha +1)}+\frac{b E_{\beta}(-x^{\beta}) t^{\alpha }}{\Gamma (\alpha +1)}+\frac{d \tau E_{\beta}(-x^{\beta}) t^{\alpha }}{\Gamma (\alpha +1)}\nonumber \\
&+& \frac{b \delta \tau E_{\beta}(-x^{\beta}) t^{2 \alpha }}{\Gamma (2 \alpha +1)}-\frac{b \delta E_{\beta}(-x^{\beta}) t^{2 \alpha }}{\Gamma (2 \alpha +1)}+\frac{b \tau E_{\beta}(-x^{\beta}) t^{2 \alpha }}{\Gamma (2 \alpha +1)}\nonumber \\ &-& \frac{b E_{\beta}(-x^{\beta}) t^{2 \alpha }}{\Gamma (2 \alpha +1)}+\frac{d \tau^2 E_{\beta}(-x^{\beta}) t^{2 \alpha }}{\Gamma (2 \alpha +1)}+\frac{b \delta \tau^2 E_{\beta}(-x^{\beta}) t^{3 \alpha }}{\Gamma (3 \alpha +1)}\nonumber \\ &-& \frac{b \delta \tau E_{\beta}(-x^{\beta}) t^{3 \alpha }}{\Gamma (3 \alpha +1)}+\frac{b \delta E_{\beta}(-x^{\beta}) t^{3 \alpha }}{\Gamma (3 \alpha +1)}+\frac{b \tau^2 E_{\beta}(-x^{\beta}) t^{3 \alpha }}{\Gamma (3 \alpha +1)}\nonumber \\ &-& \frac{b \tau E_{\beta}(-x^{\beta}) t^{3 \alpha }}{\Gamma (3 \alpha +1)}+\frac{b E_{\beta}(-x^{\beta}) t^{3 \alpha }}{\Gamma (3 \alpha +1)}+\frac{d \tau^3 E_{\beta}(-x^{\beta}) t^{3 \alpha }}{\Gamma (3 \alpha +1)}\nonumber \\ 
&+& \frac{b \delta \tau^3 E_{\beta}(-x^{\beta}) t^{4 \alpha }}{\Gamma (4 \alpha +1)}-\frac{b \delta \tau^2 E_{\beta}(-x^{\beta}) t^{4 \alpha }}{\Gamma (4 \alpha +1)}\nonumber +\frac{b \delta \tau E_{\beta}(-x^{\beta}) t^{4 \alpha }}{\Gamma (4 \alpha +1)}\nonumber \\ &-&\frac{b \delta E_{\beta}(-x^{\beta}) t^{4 \alpha }}{\Gamma (4 \alpha +1)}+ \frac{b \tau^3 E_{\beta}(-x^{\beta}) t^{4 \alpha }}{\Gamma (4 \alpha +1)}-\frac{b \tau^2 E_{\beta}(-x^{\beta}) t^{4 \alpha }}{\Gamma (4 \alpha +1)}\nonumber \\ &+&\frac{b \tau E_{\beta}(-x^{\beta}) t^{4 \alpha }}{\Gamma (4 \alpha +1)}-\frac{b E_{\beta}(-x^{\beta}) t^{4 \alpha }}{\Gamma (4 \alpha +1)}+\frac{d \tau^4 E_{\beta}(-x^{\beta}) t^{4 \alpha }}{\Gamma (4 \alpha +1)}+\cdots.
\end{eqnarray}
The compact form solution of the time and space fractional system (\ref{pde1}-\ref{pde1i}) is
\begin{eqnarray}\label{pdesol}
u_1(x,t)&=& b E_{\alpha}(-t^{\alpha})E_{\beta}(-x^{\beta}), \nonumber \\
u_2(x,t)&=& \Big(d E_{\alpha}(\tau t^{\alpha})+b \Big(\frac{1+\delta}{1+\tau}\Big)\Big[E_{\alpha}(\tau t^{\alpha})-E_{\alpha}(-t^{\alpha}) \Big]\Big)E_{\beta}(-x^{\beta}),\nonumber  \\
&&\tau \neq 1.
\end{eqnarray}
Note that for $\alpha=\beta=1$ solutions \eqref{pdesol} matches with as discussed in \cite{sahadevan2017exact}.

\begin{ex}
Consider the following time and space fractional system of PDEs
\begin{eqnarray}\label{mkpde1}
\frac{\partial^{\alpha_1}u_1}{\partial t^{\alpha_1}}&=&\frac{\partial^{3\beta}u_2^2}{\partial x^{3\beta}}+\eta \frac{\partial^{2\beta}}{\partial x^{2\beta}}\Big(u_1^2 \frac{\partial^{\beta}u_2}{\partial x^{\beta}} \Big),\nonumber \\
\frac{\partial^{\alpha_2}u_2}{\partial t^{\alpha_2}}&=&\frac{\partial^{3\beta}u_1^2}{\partial x^{3\beta}}+\zeta \frac{\partial^{2\beta}}{\partial x^{2\beta}}\Big(u_2^2 \frac{\partial^{\beta}u_1}{\partial x^{\beta}} \Big),~t>0, 0<\alpha_1,\alpha_2,\beta \leq 1,
\end{eqnarray}
where $\eta, \zeta \neq 0$ are arbitrary constants, along with the initial conditions
\begin{eqnarray}\label{mkpde1i}
u_1(x,0)= a+b x^{\beta}, ~~u_2(x,0)=c+d x^{\beta}, ~ a,b,c,d \in \mathbb{R.}
\end{eqnarray}
\end{ex}
Taking the Sumudu transform of both sides of Eqns. \eqref{mkpde1}
 \begin{eqnarray}
 S\Big[\frac{\partial^{\alpha_1}u_1}{\partial t^{\alpha_1}}\Big]&=& S\Big[\frac{\partial^{3\beta}u_2^2}{\partial x^{3\beta}}+\eta \frac{\partial^{2\beta}}{\partial x^{2\beta}}\Big(u_1^2 \frac{\partial^{\beta}u_2}{\partial x^{\beta}} \Big)\Big],\nonumber \\
 S\Big[\frac{\partial^{\alpha_2}u_2}{\partial t^{\alpha_2}}\Big]&=& S\Big[\frac{\partial^{3\beta}u_1^2}{\partial x^{3\beta}}+\zeta \frac{\partial^{2\beta}}{\partial x^{2\beta}}\Big(u_2^2 \frac{\partial^{\beta}u_1}{\partial x^{\beta}} \Big)\Big].
 \end{eqnarray}
Using the property of Sumudu transform \eqref{sproperty}, we get
\begin{eqnarray}\label{mkpde1s}
S[u_1(x,t)]&=& u_1(x,0)+\omega^{\alpha_1}S\Big[\frac{\partial^{3\beta}u_2^2}{\partial x^{3\beta}}+\eta \frac{\partial^{2\beta}}{\partial x^{2\beta}}\Big(u_1^2 \frac{\partial^{\beta}u_2}{\partial x^{\beta}} \Big)\Big],\nonumber \\
S[u_2(x,t)]&=& u_2(x,0)+\omega^{\alpha_2} S\Big[\frac{\partial^{3\beta}u_1^2}{\partial x^{3\beta}}+\zeta \frac{\partial^{2\beta}}{\partial x^{2\beta}}\Big(u_2^2 \frac{\partial^{\beta}u_1}{\partial x^{\beta}} \Big)\Big].
 \end{eqnarray}
Taking inverse Sumudu transform of both sides of Eqns. \eqref{mkpde1s}
\begin{eqnarray}
u_1(x,t)&=& S^{-1}[u_1(x,0)]+S^{-1}\Big(\omega^{\alpha_1}S\Big[\frac{\partial^{3\beta}u_2^2}{\partial x^{3\beta}}+\eta \frac{\partial^{2\beta}}{\partial x^{2\beta}}\Big(u_1^2 \frac{\partial^{\beta}u_2}{\partial x^{\beta}} \Big)\Big]\Big),\nonumber \\
u_2(x,t)&=& S^{-1}[u_2(x,0)]+S^{-1}\Big(\omega^{\alpha_2} S\Big[\frac{\partial^{3\beta}u_1^2}{\partial x^{3\beta}}+\zeta \frac{\partial^{2\beta}}{\partial x^{2\beta}}\Big(u_2^2 \frac{\partial^{\beta}u_1}{\partial x^{\beta}} \Big)\Big]\Big).\nonumber \\
\end{eqnarray}

In view of the recurrence relation \eqref{nsstim}

\begin{eqnarray}
u_{10}&=& S^{-1}[u_1(x,0)]=a+b x^{\beta}, \nonumber \\
u_{20}&=& S^{-1}[u_2(x,0)]=c+d x^{\beta}, \nonumber \\
u_{11}&=& S^{-1}\Big(\omega^{\alpha_1}S\Big[\frac{\partial^{3\beta}u_{20}^2}{\partial x^{3\beta}}+\eta \frac{\partial^{2\beta}}{\partial x^{2\beta}}\Big(u_{10}^2 \frac{\partial^{\beta}u_{20}}{\partial x^{\beta}} \Big)\Big]\Big), \nonumber \\
&=&   \frac{\eta b^2 d \Gamma (\beta +1) \Gamma (2 \beta +1) t^{\alpha _1}}{\Gamma \left(\alpha _1+1\right)}, \nonumber \\
u_{21}&=& S^{-1}\Big(\omega^{\alpha_2} S\Big[\frac{\partial^{3\beta}u_{10}^2}{\partial x^{3\beta}}+\zeta \frac{\partial^{2\beta}}{\partial x^{2\beta}}\Big(u_{20}^2 \frac{\partial^{\beta}u_{10}}{\partial x^{\beta}} \Big)\Big]\Big), \nonumber \\
&=& \frac{\zeta b d^2 \Gamma (\beta +1) \Gamma (2 \beta +1) t^{\alpha _2}}{\Gamma \left(\alpha _2+1\right)}, \nonumber \\
u_{1n}&=& 0, n\geq 2, \nonumber \\
u_{2n} &=& 0, n\geq 2.
\end{eqnarray}
Thus, the exact solution of the fractional system \eqref{mkpde1} alsong with initial conditions \eqref{mkpde1i} is given by
\begin{eqnarray}
u_1(x,t)&=& a+b x^{\beta}+\frac{ \eta b^2 d \Gamma (\beta +1) \Gamma (2 \beta +1) t^{\alpha _1}}{\Gamma \left(\alpha _1+1\right)}, \nonumber \\
u_2(x,t)&=& c+d x^{\beta}+\frac{\zeta b d^2 \Gamma (\beta +1) \Gamma (2 \beta +1) t^{\alpha _2}}{\Gamma \left(\alpha _2+1\right)}.
\end{eqnarray}

\section{Conclusions}\label{conclusion}
Sumudu transform iterative method is developed by combining Sumudu transform and DGJM \cite{daftardar2006iterative}. This approach is suitable for getting exact solutions of time and space FPDEs and as well as systems of them. We demonstrate its applicability by solving a large number of non-trivial examples. Although combinations of Sumudu transform with other decomposition methods such as HPM and ADM have been proposed in the literature\cite{singh2011homotopy,kumar2012sumudu}, the combination of Sumudu transform with DGJM gives better and more efficient method as we do not need to construct homotopy or find Adomian polynomials.

\bibliographystyle{ieeetr}
\bibliography{reference1}   
\end{document}